\documentclass[12pt]{article}%
\usepackage[hidelinks,hypertexnames=false]{hyperref}
\usepackage{amsmath,amsthm}
\usepackage{amsfonts}
\usepackage{amssymb}
\usepackage{thmtools}
\usepackage{graphicx}
\usepackage{cleveref}
\usepackage{color}%
\usepackage{enumitem}
\providecommand{\U}[1]{\protect\rule{.1in}{.1in}}
\newtheorem{theorem}{Theorem}[section]

\newtheorem{corollary}[theorem]{Corollary}

\newtheorem{lemma}[theorem]{Lemma}

\newtheorem{proposition}[theorem]{Proposition}
\newtheorem{remark}[theorem]{Remark}

\newcommand{\BIGOP}[1]{\mathop{\mathchoice{\raise-0.22em\hbox{\huge
$#1$}} {\raise-0.05em\hbox{\Large $#1$}}{\hbox{\large $#1$}}{#1}}}

\makeatletter\@addtoreset{equation}{section}\makeatother
\evensidemargin=\oddsidemargin
\newdimen\dummy
\dummy=\oddsidemargin
\crefname{subsection}{subsection}{subsections}
\renewcommand{\forall}{\text{for all }}
\newcommand{\CutC}{\mathbb{C}^{\bullet}_{\geq0}}
\newcommand{\Part}{}
\newcommand{\DtNLaplace}{\operatorname*{DtN}\left(0\right)}
\newcommand{\DtNlog}{\operatorname{DtN}^{\operatorname{log}}}
\newcommand{\SLaplace}{\operatorname{S}_{\operatorname{D}}\left(0\right)}
\newcommand{\Slog}{\operatorname{S}^{\operatorname{log}}_{\operatorname{D}}}
\newcommand{\W}[1]{{#1}}

\let\originalleft\left
\let\originalright\right
\renewcommand{\left}{\mathopen{}\mathclose\bgroup\originalleft}
\renewcommand{\right}{\aftergroup\egroup\originalright}

\begin{document}

\title{Dirichlet-to-Neumann operator for the Helmholtz problem with general wavenumbers on the $n$-sphere }
\author{Benedikt Gr\"{a}{\ss }le\thanks{Institut f\"{u}r Mathematik, Universit\"{a}t
Z\"{u}rich, Winterthurerstr.~190, CH-8057 Z\"{u}rich, Switzerland
\newline(\{benedikt.graessle, stas\}@math.uzh.ch).}
\and Stefan A. Sauter\footnotemark[1]}
\date{}
\maketitle

\begin{abstract}
This paper considers the Helmholtz problem in the exterior of a ball
with Dirichlet boundary conditions and radiation conditions imposed at
infinity. The differential Helmholtz operator depends on the complex wavenumber with
non-negative real part and is formulated for general spatial dimensions. We
prove wavenumber explicit continuity estimates of the corresponding
Dirichlet-to-Neumann ($\operatorname*{DtN}$) operator which do not deteriorate as the
complex wavenumber tends to zero.

The exterior Helmholtz problem can be equivalently reformulated on a bounded
domain with $\operatorname*{DtN}$ boundary conditions on the artificial
boundary of a ball. We derive wavenumber independent trace and Friedrichs-type 
inequalities for the solution
space in wavenumber-indexed norms. 
\end{abstract}



\section{Introduction\label{sec:Introduction}}
Acoustic scattering problems modelled by the Helmholtz equation are fundamental
for the understanding of acoustic wave propagation and the design of
efficient numerical methods. In many real-world applications 
such as design of antennas, radar detection, ultrasound imaging, and seismic wave analysis, the
computational domain is the unbounded exterior of a bounded scatterer or of an
emitting source. This unboundedness presents additional challenges in both theoretical analysis and
the design of effective solution methods.

A popular approach to
address this difficulty is to introduce a sufficiently large ball that contains the
compact, possibly, non-homogeneous material and the support of a compactly supported source term.
The boundary of this ball serves as an artificial boundary to the original exterior problem,
restricting the Helmholtz equation to its interior.
The key idea is to impose \emph{transparent boundary conditions}
on the artificial boundary such that a solution to the original
problem (restricted to the interior of the ball) coincides with the solution
on the finite domain, cf.%
~\cite{Fen:FiniteElementMethod1983,GPS:HelmholtzEquationHeterogeneous2019,GHS:StableSkeletonIntegral2025,KG:ExactNonreflectingBoundary1989,MM:FiniteElementMethod1980,Nedelec01,SW:WavenumberExplicitParametricHolomorphy2023}
to only name a few.

This approach relies on the well-posedness of the unbounded problem, which,
for exterior problems, requires appropriate \emph{conditions towards infinity},
such as the Sommerfeld radiation condition \cite{Sch:EightyYearsSommerfelds1992}. 
Under such conditions, the solution to the homogeneous Helmholtz problem in the exterior
of the ball is uniquely determined by its Dirichlet data on the artificial boundary. 
The Dirichlet-to-Neumann (DtN) operator maps a given Dirichlet trace to the 
normal derivative of the corresponding unique solution at the artificial boundary.
The transparent boundary condition --- 
also called $\operatorname*{DtN}$ boundary condition ---
on the artificial boundary
requires the Dirichlet and Neumann traces of the truncated solution
to match those of the original solution, which are expressed by the DtN operator.
In this way, the DtN operator allows one to trade a finite computational domain 
for an additional boundary condition.

The theoretical justification and practical implementation
of this approach require a thorough understanding 
of the involved DtN operator.
In particular, for the numerical discretization of
scattering problems, it is crucial to
understand the explicit behaviour with respect to
the wavenumber of the arising operators. Although
the Fredholm theory in combination with (variants
of) the unique continuation principle is a potent
tool to establish well-posedness of acoustic problems,
wavenumber explicit estimates require more refined analytical techniques.
Among them are subtle definiteness (sign) properties of the real and
imaginary part of the DtN operator%
~\cite{GPS:HelmholtzEquationHeterogeneous2019,GHS:StableSkeletonIntegral2025,MelenkSauterMathComp,SW:WavenumberExplicitParametricHolomorphy2023}.
However, these properties have only been established for purely imaginary wavenumbers in
$n=2,3$ dimensions~\cite{MonkChandlerWilde,MelenkSauterMathComp,Nedelec01},
which prevents a unified theory for general complex wavenumbers:
Existing results on Helmholtz problem based on the DtN approach 
are typically restricted to purely imaginary wavenumbers for $n\in\{2,3\}$, see, e.g.,%
~\cite{GPS:HelmholtzEquationHeterogeneous2019,SW:WavenumberExplicitParametricHolomorphy2023}
for recent developments on general coefficient Helmholtz problems.

To close this gap, this paper presents a unified mathematical
analysis of spherical DtN operators for Helmholtz problems with general wavenumbers 
in any spatial dimension and provides
a clearer mathematical understanding of the important definiteness 
characteristics of the DtN operator. 
In particular, our results enable 
a unified analysis of generalised boundary layer operators for 
Helmholtz problems with general coefficients and
wavenumbers for $n\geq2$ dimensions in~\cite{GHS:StableSkeletonIntegral2025}.
The new insights also facilitate an improved tuning of control parameters 
in numerical methods depending on the real and imaginary
part of the wavenumber such as~\cite{BLS:VariableOrderDirectional2021} in future research.%
\medskip

\noindent\textbf{Main results %
and further contributions%
.} 

\textbf{a)} 
The wavenumber-dependent DtN operator admits an explicit Fourier series
representation in terms of eigenfunctions of the
Laplace-Beltrami operator and Bessel functions.
A careful wavenumber-explicit control of the Fourier series coefficients enables sharp 
bounds on the real and imaginary part of the DtN operator for general wavenumbers 
and spatial dimensions.
Our results extend the existing theory for purely imaginary wavenumbers,
refining the bounds in~\cite[Sec.~3]{MelenkSauterMathComp} for two dimensions
and recovering the classical results \cite[Eqn.~(2.6.23)--(2.6.24)]{Nedelec01} for
three dimensions.  

\textbf{b) }The solution space of the Helmholtz equation in the truncated domain with 
transparent ($\operatorname*{DtN}$) boundary conditions on the artificial boundary 
satisfies a novel Friedrichs-type inequality in $n\geq3$ dimensions. 
This enables wavenumber-explicit Sobolev and trace estimates in the natural weighted
norms which do \emph{not} deteriorate for small complex wavenumbers with non-negative real part.
These estimates are the key to an explicit characterisation of the natural weighted trace norm for 
Helmholtz problems on exterior domains in~\cite[Thm.~3.4]{Gra:OptimalTraceNorms2025}.

\textbf{c) } The enforced condition towards infinity affects the
definition of the solution operator for the exterior Dirichlet Helmholtz problem,
whose normal derivative defines the $\operatorname*{DtN}$ operator. 
We discuss various types of radiation conditions and derive the
well-posedness of the Helmholtz problem in the exterior of a ball with
Dirichlet conditions on its boundary. This extends the results developed in
\cite{Mclean00} as well as the well-posedness results for the
$\operatorname*{DtN}$ operator in \cite{coltonkress_inverse}.

We expect that these results in the field of applied analysis
can be useful for analysing numerical discretization
methods for problems with transparent boundary conditions or approximations thereof.

\medskip\noindent
\textbf{Outline.}
\Cref{ChapPrelim} introduces the notation for the geometric setting 
and corresponding Sobolev spaces 
as well as some standard trace estimates. 

\Cref{sec:The Dirichlet-to-Neumann operator} presents the exterior Helmholtz
problem with Dirichlet boundary conditions and discusses the popular choices of the
radiation conditions towards infinity. 
The well-posedness of the resulting exterior problems is 
established by \Cref{PropProblems}. Its proof follows mainly the
theory developed in \cite{Mclean00} and is postponed to the appendix.
\Cref{thm:DtN} states the important definiteness properties of the DtN operator 
with sharp bounds on its real and imaginary parts.
The Friedrichs-type inequality described in \textbf{b)}
is formulated in \Cref{thm:C_Friedrich} and leads in \Cref{CorTrace}
to the announced wavenumber-explicit Sobolev and trace estimates. 

\Cref{sec:Spectral expansion of the DtN operator} introduces series
expansions of the $\operatorname*{DtN}$ operators for the different radiation
conditions. Sharp wavenumber-explicit estimates of the series coefficients for
non-zero wavenumbers are formulated as \Cref{thm:zmn} and proved in
\Cref{sub:Analysis of the spectral coefficients}. 
The special case of wavenumber zero (i.e.,
the Laplacian) requires extra care and is analysed in \Cref{sub:DtN_Laplace}.
This paper concludes with the proofs of \Cref{thm:DtN,thm:C_Friedrich} in \Cref{sub:Proof of DtN}.

\section{Preliminaries\label{ChapPrelim}}

The set of natural numbers is denoted by $\mathbb{N}:=\left\{  1,2,\dots
\right\}  $ and we write $\mathbb{N}_{0}:=\mathbb{N}\cup\left\{  0\right\}  $.
The following subset of the complex numbers will be needed:%
\[
\mathbb{C}_{>0}:=\left\{  z\in\mathbb{C}\mid\operatorname{Re}%
z>0\right\}  ,\quad\mathbb{C}_{\geq0}:=\left\{  z\in\mathbb{C}%
\mid\operatorname{Re}z\geq0\right\}  ,
\quad\CutC:=\mathbb{C}_{\geq0}\backslash\left\{  0\right\}  .
\]

The physical problem is formulated in $\mathbb{R}^{n}$ for general dimension
$n\in\left\{  2,3,\ldots\right\}  $. We consider Lipschitz domains
$\omega\subset\mathbb{R}^{n}$ which are open and connected (not necessarily
bounded) sets whose boundary is compact and locally the graph of a Lipschitz
function~\cite[Def.\ 3.28]{Mclean00}. The classical Lebesgue and
Sobolev spaces on $\omega$ are denoted by $L^{2}\left(  \omega\right)  $ and
$H^{\kappa}\left(  \omega\right)  $ with norms $\left\Vert \cdot\right\Vert
_{L^{2}\left(  \omega\right)  }$ and $\left\Vert \cdot\right\Vert _{H^{\kappa
}\left(  \omega\right)  }$. For $\kappa\geq0$, the space $H_{\operatorname*{loc}%
}^{\kappa}\left(  \omega\right)  $ is given by all distributions $v$ on
the compactly supported (relative to $\mathbb R^n$) and smooth functions 
$C_{\operatorname*{comp}}^{\infty}\left(  \omega\right)  $ 
such that $\varphi v\in H^{\kappa}\left(
\omega\right)  $ for all $\varphi\in C_{\operatorname*{comp}}^{\infty}\left(
\omega\right)  $. For the definition of Sobolev spaces $H^{\kappa}\left(
\W{\Gamma}\right)  $ on $\left(  n-1\right)  $-dimensional \W{Lipschitz hypersurfaces 
$\Gamma$} and their norm $\left\Vert \cdot\right\Vert _{H^{\kappa}\left(
\Gamma\right)  }$ we refer, e.g., to~\cite[pp.~96--99]{Mclean00}. 
\W{The duality-pairing that extends the bilinear version of the $L^2(\Gamma)$ scalar product is denoted by $\left\langle \bullet,\bullet\right\rangle_{\Gamma}$ and satisfies
	\begin{align*}
		\left\langle g, h\right\rangle_{\Gamma}=\int_\Gamma g\, h\;\mathrm d s
		\qquad\text{for all }g,h\in L^2(\Gamma).
	\end{align*}%
}%

We define
wavenumber-indexed norms%
\footnote{For $s=0$, the quantity $\left\Vert v\right\Vert _{H^{1}\left(
\omega\right)  ,0}=\left\Vert \nabla v\right\Vert _{L^{2}\left(
\omega\right)  }$ is only a seminorm in $H^{1}\left(  \omega\right)  $.}
 for $s\in\mathbb{C}_{\geq0}$ by
\begin{align}
\left\Vert v\right\Vert _{H^{1}\left(  \omega\right)  ,s}  &  :=\sqrt
{\left\Vert \nabla v\right\Vert _{L^{2}\left(  \omega\right)  }^{2}+\left\vert
s\right\vert ^{2}\left\Vert v\right\Vert _{L^{2}\left(  \omega\right)  }^{2}%
}
&&\text{for all }v\in H^{1}\left(  \omega\right)  ,\label{index1}\\
\left\Vert g\right\Vert _{H^{1/2}\left(  \W{\Gamma}\right)  ,s}  &
:=\sqrt{\left\Vert g\right\Vert _{H^{1/2}\left(  \W{\Gamma}\right)  }%
^{2}+\left\vert s\right\vert \left\Vert g\right\Vert _{L^{2}\left(
\W{\Gamma}\right)  }^{2}}
&&\text{for all }g\in H^{1/2}\left(
\W{\Gamma}\right)  . \label{index2}%
\end{align}
An immediate consequence is the estimate $\Vert\bullet\Vert_{H^{1/2}\left(
\W{\Gamma}\right)  }\leq\Vert\bullet\Vert_{H^{1/2}\left(  \W{\Gamma}\right)  ,s}$.
The space 
$H^{1}\left(  \omega,\Delta\right)  :=\left\{  u\in H^{1}\left(  \omega\right)
\mid\Delta u\in L^{2}\left(  \omega\right)  \right\}$ is equipped with the (semi-)norms
\begin{align*}
	\left\Vert {v}\right\Vert _{H^{1}\left(
			\omega,\Delta\right)  }&:=
			\sqrt{\left\Vert {v}\right\Vert _{H^{1}\left(
					\omega\right)  }^{2}+\left\Vert \Delta{v}\right\Vert _{L^{2}\left(
				\omega\right)  }^{2}}
								   &&\text{for all }v\in H^1(\omega,\Delta),\\
\left\vert {v}\right\vert
_{H^{1}\left(  \omega,\Delta\right)  }&:=\sqrt{\left\Vert \nabla{v}
\right\Vert _{L^{2}\left(  \omega\right)  }^{2}+\left\Vert \Delta
{v}\right\Vert _{L^{2}\left(  \omega\right)  }^{2}}
									  &&\text{for all }v\in H^1(\omega,\Delta).
\end{align*}
It is well-known
that the outer unit normal $n(x)\in \mathbb{S}_1$ of the Lipschitz domain $\omega\subset\mathbb{R}^n$
exists for almost every $x\in\partial\omega$.
The classical normal derivative $\partial/\partial n:\overline{\omega}
\rightarrow \Gamma$ 
on some relatively open connectivity component $\Gamma\subset\partial\omega$
of the boundary can be continuously extended to a
continuous mapping $\gamma_{1}:H^{1}\left(  \omega,\Delta\right)  \rightarrow
H^{-1/2}\left(  \Gamma\right)  $ with
\begin{equation}
\left\Vert \gamma_{1}v\right\Vert _{H^{-1/2}\left(  \Gamma\right)  }\leq
C_{\omega,1}\left\vert v\right\vert _{H^{1}\left(  \omega,\Delta\right)
}\quad\forall v\in H^{1}\left(  \omega,\Delta\right)  .\label{normaltraceest}%
\end{equation}
Also the following standard and multiplicative trace estimates are well known:%
\begin{equation}
\left\Vert v\right\Vert _{H^{1/2}\left(  \Gamma\right)  }\leq C_{\omega
,0}\left\Vert v\right\Vert _{H^{1}\left(  \omega\right)  }\quad\text{and\quad
}\left\Vert v\right\Vert _{L^{2}\left(  \Gamma\right)  }\leq C_{\omega
,0}\left\Vert v\right\Vert _{L^{2}\left(  \omega\right)  }^{1/2}\left\Vert
v\right\Vert _{H^{1}\left(  \omega\right)  }^{1/2}\quad\forall v\in
H^{1}\left(  \omega\right)  .\label{standtrace}%
\end{equation}
\W{
	We define the Hilbert space $H^{1}\left(  \omega,\Gamma\right)$ as the completion of $H^{1}\left(  \omega,\Delta\right)  $ with respect to the norm%
\[
\left\Vert \bullet\right\Vert _{H^1(\omega,\Gamma)}:=\sqrt{\left\Vert \bullet\right\Vert
_{H^{1}\left(  \omega\right)  }^{2}+\left\Vert \gamma_{1}\bullet\right\Vert
_{H^{-1/2}(\Gamma)}^{2}}
\]
and note that $\gamma_1:H^1(\omega,\Gamma)\to H^{-1/2}(\Gamma)$ is continuous by construction.}
The space $L(H_1; H_2)$ of continuous linear maps $A:H_1\to H_2$ between Hilbert spaces 
$H_1$ and $H_2$ is equipped with the operator norm topology.
\W{An operator sequence $\{A_n\}_{n\in\mathbb N}\subset L(H_1;H_2)$ converges strongly to $A\in L(H_1;H_2)$ as $n\to \infty$, if and only if $A_n x\to Ax$ converges in the $H_2$-norm for every $x\in H_1$.
A map $A:D\subset\mathbb C\to L(H_1;H_2)$ is strongly continuous, if and only if $A(s)$ converges strongly to $A(s_0)$ for all $s,s_0\in D$ as $s\to s_0$.
For the Fr\'echet space $H^1_{\mathrm{loc}}(\omega)$,
}%
we set%
\begin{align*}
	L(H_1; H^1_{\mathrm{loc}}(\omega))
	 :=\{ A:H_1\to H^1_{\mathrm{loc}}(\omega)\ 
		 \vert\ \varphi A\in L(H_1; H^1(\omega))\text{ for all }\varphi\in
	C^\infty_{\mathrm{comp}}(\omega)\}.%
\end{align*}

Let $B_{R}\subset\mathbb{R}^{n}$ denote the open ball of radius $R>0$ at the origin. 
The outer unit normal vector $x/\left\Vert x\right\Vert$ 
on the boundary $S_{R}=\partial B_R$ points into the unbounded complement 
$B_R^+:=\mathbb{R}^n\setminus\overline{B_R}$ of $B_{R}$ and its derivative in this direction 
is denoted by $\partial_{r}$. 
The notation $|\bullet|$ is context-sensitive and may refer to the Lebesgue measure $|\omega|$ of a
bounded measurable $n$-dimensional set $\omega\subset\mathbb{R}^{n}$, 
the surface measure $|\Gamma|$ of an $\left(
n-1\right)  $-dimensional manifold $\Gamma\subset\mathbb{R}^{n}$, and 
the cardinality $|J|$ of a countable set $J$.

\begin{remark}
\label{RemNormalHharm}Let $s\in\mathbb{C}_{\geq0}$. Note that any function
$u\in H_{\operatorname*{loc}}^{1}\left(  B_{R}^{+}\right)  $ with $\Delta
u+s^{2}u=0$ in a bounded neighborhood $\omega$ of $S_{R}$ in $B_{R}^{+}$ satisfies
$\left.  u\right\vert _{\omega}\in H^{1}\left(  \omega,\Delta\right)  $ and
the normal derivative $\partial_{r}u$ of $u$ on $S_{R}$ is well defined in
$H^{-1/2}\left(  S_{R}\right)  $.
\end{remark}

\section{The Dirichlet-to-Neumann operator on a
sphere\label{sec:The Dirichlet-to-Neumann operator}}

This section discusses the well-posedness of the Dirichlet Helmholtz problem on the
exterior of the ball $B_R$
and the properties of the associated
Dirichlet-to-Neumann operator (DtN) on the
sphere $S_{R}$.

\subsection{The Dirichlet Helmholtz problem on the exterior of a ball\label{sub:The DtN operator for the Helmholtz problem on a sphere}}
Let $s\in\mathbb{C}_{\geq0}$ denote the wavenumber
and consider the exterior Dirichlet Helmholtz problem %
\begin{equation}
\text{find }u\in H_{\operatorname*{loc}}^{1}\left(  B_{R}^{+}\right)
\quad\text{such that\quad}\left\{
\begin{array}
[c]{rl}%
-\Delta u+s^{2}u=0 & \text{in }B_{R}^{+}= \mathbb{R}\setminus\overline{B_R},\\
u=g & \text{on }S_{R}=\partial B_R.
\end{array}
\right.  \label{eqn:exterior_problem}%
\end{equation}
These equations can be closed by prescribing the behaviour of $u$ by some
\emph{conditions towards infinity}
(denoted by $\operatorname*{IC}_{s}$). These
conditions are by no means unique; the first version is given by (\ref{SRCa}),
(\ref{SRCb}) and denoted by $\operatorname*{IC}_{s}^{\operatorname*{strong}}$%
\begin{subequations}
\label{SRC}
\end{subequations}%
\begin{equation}
\left\{
\begin{array}
[c]{rl}%
\left\vert u\left(  r\xi\right)  \right\vert  & \leq Cr^{\mu}\\
\left\vert \frac{\partial u}{\partial r}\left(  r\xi\right)  +su\left(
r\xi\right)  \right\vert  & \leq Cr^{\mu-1}%
\end{array}
\right\}  \quad\text{as }r\rightarrow\infty\text{ uniformly in }\xi\in
S_{1}.\tag{%
\ref{SRC}%
a}\label{SRCa}%
\end{equation}
The rate $\mu$ depends on the spatial dimension $n$ and the wavenumber $s:$%
\begin{equation}
\mu:=\left\{
\begin{array}
[c]{ll}%
\frac{1-n}{2} & \text{if }s\in\CutC,\\
2-n & \text{if }s=0.
\end{array}
\right.  \tag{%
\ref{SRC}%
b}\label{SRCb}%
\end{equation}

\begin{remark}[infinity conditions]
\label{RemICs}There exist alternative conditions towards infinity in the
literature\footnote{The superscript \textquotedblleft$\operatorname*{Green}%
$\textquotedblright\ indicates that for $\operatorname*{IC}_{s}%
^{\operatorname*{Green}}$ a Green's representation theorem holds \textit{under
certain assumptions} while the superscript \textquotedblleft%
$\operatorname*{isom}$\textquotedblright\ indicates that these conditions
\textit{always} lead to well-posed problems.}.

\begin{enumerate}
\item Condition $\operatorname*{IC}_{s}^{\operatorname*{Green}}$ is 

\begin{enumerate}
\item for $s\in\CutC$ given by%
\begin{equation}
\lim_{r\rightarrow\infty}r^{(n-1)/2}\left(  \partial_{r}{u}\left(
r\xi\right)  +s{u}\left(  r\xi\right)  \right)  =0\quad\text{uniformly in
}\xi\in S_{1},\label{SRCc}%
\end{equation}

\item for $s=0$ and

\begin{enumerate}
\item $n\geq3$ given by%
\[
\left\vert u\left(  r\xi\right)  \right\vert \leq Cr^{2-n}\quad\text{as
}r\rightarrow\infty\text{ uniformly in }\xi\in S_{1},
\]

\item $n=2$ given by%
\begin{equation}
\exists b\in\mathbb{R}\qquad u\left(  r\xi\right)  -b\log r\leq Cr^{-1}%
\quad\text{as }r\rightarrow\infty\text{ uniformly in }\xi\in S_{1}.
\label{SRCd}%
\end{equation}

\end{enumerate}
\end{enumerate}

\item Condition $\operatorname*{IC}_{s}^{\operatorname*{isom}}$ is the same as
condition $\operatorname*{IC}_{s}^{\operatorname*{Green}}$ for $s\in
\mathbb{C}_{\geq0}$ and $\left(  s,n\right)  \neq\left(  0,2\right)  $ while
for $s=0$ and $n=2$, condition $\operatorname*{IC}_{0}^{\operatorname*{isom}}$
is given by%
\[
\left\vert u\left(  r\xi\right)  \right\vert \leq C\quad\text{as
}r\rightarrow\infty\text{ uniformly in }\xi\subset S_{1}.
\]

\end{enumerate}
\end{remark}

In the following proposition we state the well-posedness of these exterior
Dirichlet problems.
We say that problem~\eqref{eqn:exterior_problem}
is well posed for $\operatorname*{IC}_{s}$, if the corresponding solution operator 
$\operatorname{S_D}\in L(H^{1/2}(S_R); H^1_{\mathrm{loc}}(B_R^+))$
that assigns each Dirichlet data $g\in H^{1/2}(S_R)$ the unique solution 
$\operatorname{S_D}g\in H^{1}_{\mathrm{loc}}(B_R^+)$ with $\operatorname{IC}_s$
to~\eqref{eqn:exterior_problem}, is
well-defined and continuous.

\begin{proposition}[well-posedness]
\label{PropProblems}\qquad

\begin{enumerate}
\item Let $s\in\CutC$.

\begin{enumerate}
\item Problem (\ref{eqn:exterior_problem}) with $\operatorname*{IC}%
_{s}^{\operatorname*{Green}}$ is well posed with a continuous
solution operator denoted by $\operatorname*{S}_{\operatorname*{D}}\left(  s\right)
\in L(H^{1/2}\left(  S_{R}\right); H_{\operatorname*{loc}}^{1}%
(B_{R}^{+}))$.

\item For any $g\in H^{1/2}\left(  \Gamma\right)  $ the solution
$\operatorname*{S}_{\operatorname*{D}}\left(  s\right)  g$ also satisfies the
stronger condition $\operatorname*{IC}_{s}^{\operatorname*{strong}}$. 
Hence (\ref{eqn:exterior_problem}) with $\operatorname*{IC}%
_{s}^{\operatorname*{strong}}$ is well posed with the solution operator
$\operatorname*{S}_{\operatorname*{D}}\left(  s\right)$.
\end{enumerate}

\item Let $s=0$.

\begin{enumerate}
\item Then problem (\ref{eqn:exterior_problem}) with $\operatorname*{IC}%
_{0}^{\operatorname*{strong}}$ is well posed with a continuous
solution operator denoted by ${\SLaplace}%
\in L(H^{1/2}\left(  S_{R}\right); H_{\operatorname*{loc}}^{1}
(B_{R}^{+}))$.

\item Let $s=0$ and $n=2$. If and only if $R\neq1$, the problem
(\ref{eqn:exterior_problem}) with $\operatorname*{IC}_{0}%
^{\operatorname*{Green}}$ is well posed with a continuous solution operator denoted by ${\Slog}
\in L(H^{1/2}\left(  S_{R}\right);  H_{\operatorname*{loc}}^{1}%
(B_{R}^{+}))$.
\end{enumerate}
\end{enumerate}
\end{proposition}

Major part of the proof of \Cref{PropProblems} given in the appendix 
is a combination of results which are worked out in \cite{Mclean00} in more generality.
\begin{remark}[Lipschitz domains]
	The results from \Cref{PropProblems} generalise to Dirichlet Helmholtz problems
	on the exterior of general bounded Lipschitz domains with the results in%
	~\cite{Mclean00}.
	In this case, the condition $R\ne1$ for the Laplace case in $n=2$ dimensions
	in \Cref{PropProblems}.2b has to be replaced by 
	the condition that the \emph{logarithmic capacity}
	$\operatorname{Cap}_{\Gamma}$ defined in~\cite[p.~264]{Mclean00}
	does \emph{not} equal one, i.e., $\operatorname{Cap}_{\Gamma}\ne1$.
\end{remark}

\subsection{The DtN operator for the Helmholtz problem on a
sphere}%
\label{sub:DtN_Helmholtz}

The solution operators for the exterior Dirichlet Helmholtz problems 
from \Cref{PropProblems} define the \emph{Dirichlet-to-Neumann} operator. 
Recall that the
normal derivative $\partial_{r}$ at $S_{R}$ is well defined for any solution 
to~(\ref{eqn:exterior_problem}) by \Cref{RemNormalHharm}. 
For any $s\in\mathbb{C}_{\geq0}$, the operator $\operatorname*{DtN}%
\left(  s\right)  :H^{1/2}\left(  S_{R}\right)  \rightarrow H^{-1/2}\left(
S_{R}\right)  $ is given by%
\[
\operatorname*{DtN}\left(  s\right)  g:=\partial_{r}\operatorname*{S}%
\nolimits_{\operatorname*{D}}\left(  s\right)  g\qquad\text{for all }g\in
H^{1/2}\left(  S_{R}\right)  .
\]
In the special case $s=0$ and $R\ne 1$ in $n=2$ dimensions, we additionally set%
\[
\DtNlog  g:=\partial_{r}{\Slog}  g
  \qquad\text{for all }g\in H^{1/2}\left(  S_{R}\right)  .
\]
For functions $v\in H_{\operatorname*{loc}}^{1}\left(  B_{R}^{+}\right)$, we
write $\operatorname*{DtN}\left(  s\right)  v$ for $\operatorname*{DtN}\left(
s\right)  \left(  \left.  v\right\vert _{S_{R}}\right)  $ and analogously for
${\DtNLaplace}$.

The Dirichlet-to-Neumann operators allows us to characterize the solutions to
the exterior problems stated in \Cref{PropProblems} by its traces
on the sphere $S_{R}$:
For $s\in\mathbb{C}_{\geq0}$, their definition implies for any Helmholtz-harmonic
function $u\in H_{\operatorname*{loc}}^{1}(B_{R}^{+})$ (i.e., $-\Delta
u+s^{2}u=0$ in $B_{R}^{+}$) that

\begin{itemize}
\item condition $\operatorname*{IC}_{s}^{\operatorname*{strong}}$ holds if and only
if its normal derivative satisfies
\[%
\begin{array}
[c]{cl}%
\partial_{r}u=\operatorname*{DtN}\left(  s\right)    u   &
\text{for }s\in\mathbb{C}_{\geq0}\text{ and }u=\operatorname{S}_{\operatorname{D}}\left(s\right)(u|_{S_R}),
\end{array}
\]
\item condition $\operatorname*{IC}_{0}^{\operatorname*{Green}}$ holds for
$n=2, s=0$, and $R\neq1$ if and only if%
\[
\partial_{r}u=\DtNlog  \left(  u|_{S_{R}}\right)
\quad\text{for }u={\Slog}  \left(  u|_{S_{R}}\right)  .
\]

\end{itemize}

The first main result states the key properties of the DtN operator 
including bounds on its real and imaginary parts\W{, and continuity in the strong operator topology}.
The symmetry and the boundedness, 
given here in terms of the natural wavenumber-indexed norms,
are more standard and included with elementary proofs for completeness.
\begin{theorem} [properties of $\operatorname*{DtN}$]\label{thm:DtN}The 
operator 
$\operatorname*{DtN}\left(  s\right)  \in L(H^{1/2}\left(  S_{R}\right);H^{-1/2}\left(  S_{R}\right))$
\W{is strongly continuous} on the wavenumber $s\in\CutC$
and admits a \W{strongly} continuous extension to $\mathbb{C}_{\geq0}$,%
\W{
\begin{equation}
{\DtNLaplace}g:= \lim_{\substack{s\in\CutC\\s\rightarrow
0}}\operatorname*{DtN}\left(  s\right)g
\qquad\text{for all }g\in H^{1/2}(S_R).\label{limtrel}%
\end{equation}%
}%
Any $g,h\in H^{1/2}\left(  S_{R}\right)  $ satisfy

\begin{enumerate}[label=(\roman*)]
\item
$\displaystyle
\left\langle \operatorname*{DtN}\left(  s\right)  g,\overline{h}%
\right\rangle _{S_{R}}=\left\langle g,\operatorname*{DtN}\left(  s\right)  \overline{h}%
\right\rangle _{S_{R}}=\left\langle g,\overline{\operatorname*{DtN}\left(
\overline{s}\right)  h}\right\rangle _{S_{R}},$
\item $\displaystyle\left\vert \left\langle \operatorname*{DtN}\left(  s\right)
g,\overline{h}\right\rangle _{S_{R}}\right\vert \leq C\left(  n\right)
\left\Vert g\right\Vert _{H^{1/2}\left(  S_{R}\right)  ,s}\left\Vert
h\right\Vert _{H^{1/2}\left(  S_{R}\right)  ,s},$ 

\item $\displaystyle\frac{n-2}{2}\left\Vert g\right\Vert _{L^{2}\left(  S_{R}\right)  }^{2}%
\leq-\operatorname{Re}\left(  \left\langle \operatorname*{DtN}\left(
s\right)  g,\overline{g}\right\rangle _{S_{R}}\right)  \leq C\left(  n\right)
\left\Vert g\right\Vert _{H^{1/2}\left(  S_{R}\right)  ,s}^{2},$

\item $\displaystyle0<\mp\operatorname{Im}\left(  \left\langle \operatorname*{DtN}\left(
s\right)  g,\overline{g}\right\rangle _{S_{R}}\right)  \leq C(n)\Vert
g\Vert_{H^{1/2}\left(  S_{R}\right)  ,s}^{2}\hfill\text{whenever }%
\pm\operatorname{Im}s>0\text{ and }g\neq0.$

\end{enumerate}
The constant $C(n)>0$ exclusively depends on $R$ and the dimension $n\geq2$.
\end{theorem}
The second main result is a
Friedrichs-type inequality for functions with homogeneous DtN boundary
conditions on $S_{R}$ in $n\geq 3$ dimensions.
The proofs of \Cref{thm:DtN,thm:C_Friedrich}
follow in \Cref{sub:Proof of DtN} below based
an explicit representation of
the $\operatorname*{DtN}\left(  s\right)$ operator derived in
\Cref{sec:Spectral expansion of the DtN operator} in terms of its Fourier coefficients.

\begin{theorem}
[Friedrichs-type inequality]\label{thm:C_Friedrich} Consider a bounded
Lipschitz domain $\Omega\subset\mathbb{R}^{n}$ in $n\geq3$ dimensions such
that $S_{R}\subset\partial\Omega$ is relatively closed in $\partial\Omega$.
Any $v\in H^{1}\left(  \Omega,S_{R}\right)  $ with
$\partial_{r}v=\operatorname*{DtN}\left(  s\right)  v$ on $S_{R}$ for some
$s\in\mathbb{C}_{\geq0}$ satisfies%
\begin{equation}
\left\Vert v\right\Vert _{L^{2}\left(  \Omega\right)  }\leq
C_{\operatorname*{F}}\sqrt{\Vert\nabla v\Vert_{L^{2}\left(  \Omega\right)
}^{2}+\left\Vert \partial_{r}v\right\Vert _{H^{-1/2}\left(  S_{R}\right)
}^{2}}.\label{eqn:C_F}%
\end{equation}
The constant $C_{\operatorname*{F}}>0$ is independent of $s$ and exclusively
depends on $\Omega$ and $R$.
\end{theorem}%
\W{We remark that in $n=2$ dimensions, the constants $\mathbb{R}\subset \mathrm{ker}\operatorname{DtN}(0)$ lie in the kernel of the DtN operator at $s=0$ and prevent an estimate of the form~\eqref{eqn:C_F} uniformly in $s\in\mathbb{C}_{\geq0}$.} 
A corollary for Helmholtz-harmonic functions on a Lipschitz domain $\Omega$
concludes this section.
Recall the constants $C_{\omega,0}, C_{\omega,1}$, and $C_{\mathrm{F}}$
from~\eqref{normaltraceest}--\eqref{standtrace}
and \Cref{thm:C_Friedrich}.

\begin{corollary}
\label{CorTrace}Consider a bounded Lipschitz domain $\Omega\subset\mathbb{R}^{n}$,
$n\geq3$ and let $C_{\operatorname*{F}}$ be as in \Cref{thm:C_Friedrich}. 
Any $v\in H^{1}\left(  \Omega\right)  $ with the
properties%
\[%
\begin{array}
[c]{cl}%
\left.  \left(  -\Delta v+s^{2}v\right)  \right\vert _{\omega}=0 & \text{for
some }\omega\subset\Omega\text{ such that }S_{R}\text{ is a relatively closed
subset of }\partial\omega,\\
\partial_{r}v=\operatorname*{DtN}\left(  s\right)    v   &
\text{on }S_{R}\text{ for some }s\in\mathbb{C}_{\geq0}%
\end{array}
\]
satisfies%
\begin{align}
\left\Vert v\right\Vert _{L^{2}\left(  \Omega\right)  } &  \leq\max\left\{
1,\sqrt{1+C_{\omega,1}^{2}}C_{\operatorname*{F}}\right\}  \left\Vert
v\right\Vert _{H^{1}\left(  \Omega\right)  ,s},\label{eqn:C_F_s}\\
\left\Vert v\right\Vert _{H^{1/2}\left(  S_{R}\right)  ,s} &  \leq
C_{\operatorname*{tr}}\left\Vert v\right\Vert _{H^{1}\left(  \Omega\right)
,s}.\label{esttr}%
\end{align}
In particular, if
$\left\vert s\right\vert \leq1/\sqrt{2\max\{1,C_{\omega,1}C_{\operatorname*{F}%
}\}}$, then
\begin{align}\label{eqn:C_F_gradient}
\left\Vert v\right\Vert _{L^{2}\left(  \Omega\right)  }\leq
2C_{\operatorname{F}}\sqrt{1+C_{\omega,1}^{2}}
\left\Vert \nabla v\right\Vert _{L^{2}\left(  \Omega\right)  } .
\end{align}
The constant $C_{\operatorname*{tr}}>0$ exclusively depends on the trace
constants $C_{\Omega,0}$, $C_{\omega,1}$ as well as on $C_{\operatorname*{F}}$.
\end{corollary}

\begin{proof}
	The estimate \W{\eqref{eqn:C_F_s}} is trivial for $\left\vert s\right\vert \geq1$
by definition of the norm $\left\Vert \bullet\right\Vert _{H^{1}\left(
\Omega\right)  ,s}$.
By assumption, the connected sphere $S_R\subset\partial\omega$ 
is relatively closed in $\partial\omega$
and thereby with a positive distance from $\partial\omega\setminus S_R$.
For $\left\vert s\right\vert \leq1$,
\Cref{thm:C_Friedrich} and~\eqref{normaltraceest} provide
\begin{align}
\left\Vert v\right\Vert _{L^{2}\left(  \Omega\right)  } &  
\leq C_{\operatorname*{F}}\sqrt{\left(  1+C_{\omega,1}%
^{2}\right)  \left\Vert \nabla v\right\Vert _{L^{2}\left(  \Omega\right)
}^{2}+C_{\omega,1}^{2}\left\Vert \Delta v\right\Vert _{L^{2}\left(
\omega\right)  }^{2}}.\label{vL2esta}
\end{align}
Since $\Delta v=s^2 v$ in $\omega$ and $|s|^4\leq |s|^2$,~\W{\eqref{eqn:C_F_s}}
follows for $|s|\leq1$ 
from
\begin{align}
\left\Vert v\right\Vert _{L^{2}\left(  \Omega\right)  } &  
  \leq C_{\operatorname*{F}}\sqrt{\left(  1+C_{\omega,1}^{2}\right)
\left\Vert \nabla v\right\Vert _{L^{2}\left(  \Omega\right)  }^{2}%
+C_{\omega,1}^{2}\left\vert s\right\vert ^{2}\left\Vert v\right\Vert
_{L^{2}\left(  \omega\right)  }^{2}}\leq C_{\operatorname*{F}}\sqrt
{1+C_{\omega,1}^{2}}\left\Vert v\right\Vert _{H^{1}\left(  \Omega\right)
,s}.\label{vL2estb}%
\end{align}
Let $\max\left\{  1,C_{\omega,1}C_{\operatorname*{F}}\right\}
\left\vert s\right\vert ^{2}\leq1/2$ (which implies $\left\vert s\right\vert
\leq1$) and observe from $\Delta v=s^{2}v$ in $\omega$ that%
\begin{equation}
\left\Vert \Delta v\right\Vert _{L^{2}\left(  \omega\right)  }=\left\vert
s\right\vert ^{2}\left\Vert v\right\Vert _{L^{2}(\omega)}\leq|s|^{2}\left\Vert
v\right\Vert _{L^{2}\left(  \Omega\right)  }.\label{deltavfb}%
\end{equation}
This and~\eqref{vL2esta} establish the estimate~\eqref{eqn:C_F_gradient} of the 
$L^{2}\left(  \Omega\right)  $ norm by
\begin{align*}
\left\Vert v\right\Vert _{L^{2}\left(  \Omega\right)  }/2 &  \leq\left(
1-\left\vert s\right\vert ^{2}\max\left\{  1,C_{\omega,1}C_{\operatorname*{F}%
}\right\}  \right)  \left\Vert v\right\Vert _{L^{2}\left(  \Omega\right)
}\overset{\text{(\ref{deltavfb})}}{\leq}\left\Vert v\right\Vert _{L^{2}\left(
\Omega\right)  }-\max\left\{  1,C_{\omega,1}C_{\operatorname*{F}}\right\}
\left\Vert \Delta v\right\Vert _{L^{2}\left(  \omega\right)  }\\
&  \leq C_{\operatorname*{F}}
\sqrt{1+C_{\omega,1}^{2}}\left\Vert \nabla v\right\Vert _{L^{2}\left(
\Omega\right)  }.
\end{align*}
It remains to prove the trace estimate in (\ref{esttr}). An immediate
consequence (known from~\cite[Cor.~3.2]{MelenkSauterMathComp})
of a Young inequality and the trace estimate~\eqref{standtrace} reveals
\[
|s|\left\Vert v\right\Vert _{L^{2}\left(  S_{R}\right)  }^{2}\leq
\frac{C_{\Omega,0}^{2}}{2}\left(  \left\Vert v\right\Vert _{H^{1}\left(
\Omega\right)  }^{2}+|s|^{2}\left\Vert v\right\Vert _{L^{2}\left(
\Omega\right)  }^{2}\right)  .
\]
Hence (\ref{esttr}) \W{follows from~\eqref{eqn:C_F_s}} with $C_{\operatorname*{tr}}^{2}=3\left(1+\max
\left\{1,\left(  1+C_{\omega,1}^{2}\right)  C_{\operatorname*{F}}^{2}%
\right\}\right)C_{\Omega,0}^{2}/2$ \W{as}%
\[
\left\Vert v\right\Vert _{H^{1/2}\left(  S_{R}\right)  ,s}^{2}\leq
\frac{3C_{\Omega,0}^{2}}{2}\left\Vert v\right\Vert _{H^{1}\left(
\Omega\right)  }^{2}+\frac{C_{\Omega,0}^{2}}{2}|s|^{2}\left\Vert v\right\Vert
_{L^{2}\left(  \Omega\right)  }^{2}\leq C_{\operatorname*{tr}}^{2}\left\Vert
v\right\Vert _{H^{1}\left(  \Omega\right)  ,s}^{2}.\qedhere
\]%
\end{proof}

\begin{lemma}
[special case for $n=2$]\label{lem:two dimensions}Let $s=0$, $n=2$, and $R>1$. The
difference of ${\DtNLaplace}$ for problem
(\ref{eqn:exterior_problem})--(\ref{SRC}) and
$\DtNlog  $ for problem
(\ref{eqn:exterior_problem}) and (\ref{SRCd}) is constant, namely%
\[
{\DtNLaplace}g-\DtNlog
g=-\frac{1}{R\log R}\left(  \frac{1}{\left\vert S_{R}\right\vert }\int_{S_{R}%
}g\right)  \qquad\text{for all }g\in H^{1/2}\left(  S_{R}\right)  .
\]

\end{lemma}

The proof of \Cref{lem:two dimensions} follows from an explicit Fourier
expansion of the solutions ${\SLaplace}g$
and ${\Slog} g $ in \Cref{sub:DtN_Laplace}.

\section{Spectral expansion of the DtN operator}
\label{sec:Spectral expansion of the DtN operator}
The $\operatorname*{DtN}$ operator has an
explicit series expansion in $n$-dimensional polar coordinates in terms of 
Laplace-Beltrami eigenfunctions known for $s=\operatorname{i}\mathbb R$, e.g., from
\cite{MonkChandlerWilde} for $n=2$ and 
\cite[Eqn.~(2.6.92)]{Nedelec01} for $n=3$ dimensions.
An extension to general wavenumbers $s\in\mathbb{C}_{\geq0}$ and dimensions $n\geq2$ can be derived from~\cite[Chap.~9]{Mclean00}. 

\subsection{Spectral representation of the Helmholtz DtN}
\label{sub:Spectral expansion of the Helmholtz DtN operator}
The spherical part of the DtN operator is expanded by an orthonormal basis
of the eigenfunctions of the Laplace-Beltrami operator, whose properties are briefly 
recalled, e.g., from~\cite[Sec.\ 22.2]{Shubin2001}.
The Laplace-Beltrami operator $-\Delta_{S_{1}}$ on the unit sphere
$S_{1}\subset\mathbb{R}^{n}$ is defined by%
\[
-\Delta_{S_{1}}g:=-\left.  \left(  \Delta G\right)  \right\vert _{S_{1}}%
\]
with the extension of $g$ to $\mathbb{R}^{n}\backslash\left\{  0\right\}  $ by
$G\left(  x\right)  :=g\left(  x/\left\Vert x\right\Vert \right)  $. The
eigenfunctions $Y_{m,j}$ (also called spherical harmonics) of the Laplace-Beltrami operator are given by
\begin{equation}
-\Delta_{S_{1}}Y_{m,j}=\lambda_{m}Y_{m,j},\quad m\in\mathbb{N}_{0}\text{,
}j\in J_{m}, \label{defiotam}%
\end{equation}
for some finite index set $J_{m}$ which corresponds to the multiplicity of the
eigenvalue $\lambda_{m}$ for $m\in\mathbb{N}_{0}:=\mathbb{N}\cup\{0\}$. We
assume an ascending enumeration of the eigenvalues $\lambda_{m}$, i.e.,
\[
0=\lambda_{0}<\lambda_{1}<\ldots\overset{m\rightarrow\infty}{\longrightarrow
}\infty,
\]
and that the eigenfunctions are chosen such that $\{Y_{m,j}\}_{j\in J_{m}%
,m\in\mathbb{N}_{0}}$ forms an orthonormal basis of $L^{2}\left(  S_{1}\right)
$. Explicitly it holds~\cite[Thm.~22.1]{Shubin2001} that
\begin{equation}
\lambda_{m}=m\left(  m+n-2\right)  \quad\text{with multiplicity }|J_{m}%
|=\frac{2m+n-2}{m+n-2}\binom{m+n-2}{n-2}. \label{defiotamnumber}%
\end{equation}

Next, we rewrite the exterior Helmholtz problem (\ref{eqn:exterior_problem})
in spherical coordinates $x=r\xi$, where $r:=\left\Vert x\right\Vert $ and
$\xi:=x/r\in S_{1}$. Let $x=\psi\left(  r,\xi\right)  $ denote the associated
transformation and write $\widehat{w}=w\circ\psi$ for a function $w$ in
Cartesian coordinates. The Laplace operator in spherical coordinates takes the
form
\[
\widehat{\Delta v}=\left(  \Delta v\right)  \circ\psi=\frac{1}{r^{n-1}%
}\partial_{r}\left(  r^{n-1}\partial_{r}\widehat{v}\right)  +\frac{1}{r^{2}}%
\Delta_{S_{1}}\widehat{v}.
\]
By expanding the solution $\widehat{u}=u\circ\psi$ to (\ref{eqn:exterior_problem})
by a Fourier series in terms of spherical harmonics
\begin{equation}
\widehat{u}\left(  r,\xi\right)  =\sum_{m\in\mathbb{N}_{0}}\sum_{j\in J_{m}}%
\widehat{u}_{m,j}\left(  r\right)  Y_{m,j}\left(  \xi\right)  ,\label{uhutfourier}%
\end{equation}
the Helmholtz equation (\ref{eqn:exterior_problem}) translates to an ODE for
the Fourier coefficients $\widehat{u}_{m,j}\left(  r\right)  $:%
\begin{equation}
-\frac{1}{r^{n-1}}\partial_{r}\left(  r^{n-1}\partial_{r}\widehat{u}_{m,j}\left(
r\right)  \right)  +\left(  \frac{m\left(  m+n-2\right)  }{r^{2}}%
+s^{2}\right)  \widehat{u}_{m,j}\left(  r\right)  =0\quad\forall
r>R.\label{locPDE}%
\end{equation}
Let $I_{\mu}$ and $K_{\mu}$ denote the standard modified Bessel functions of
order $\mu\in\mathbb{R}$ (e.g., from~\cite[\S10]{NIST:DLMF}) and introduce
the dimension-dependent parameter%
\begin{equation}
\nu=\frac{n-2}{2}\in\frac{1}{2}\mathbb{N}_{0}.\label{eqn:nu_def}%
\end{equation}
The ODE~(\ref{locPDE}) is closely related to the (modified) Bessel's equation
(see, e.g.,~\cite[p.~278]{Mclean00} for details).
For any $s\in\CutC$, a fundamental system is
given  by $f_{m,\nu}%
^{(1)}\left(  sr\right)  $ and $f_{m,\nu}^{(2)}\left(  sr\right)  $ with%
\begin{equation}
f_{m,\nu}^{(1)}\left(  z\right)  :=\sqrt{\frac{\pi}{2}}\frac{K_{_{m+\nu}%
}\left(  z\right)  }{z^{\nu}}\qquad\text{and}\qquad f_{m,\nu}^{(2)}\left(
z\right)  :=\frac{I_{_{m+\nu}}\left(  z\right)  }{\sqrt{2\pi}z^{\nu}}%
\qquad\text{for all }z\in\CutC.\label{deffm12}%
\end{equation}
Hence there exists coefficients $A_{m,j,\nu},B_{m,j,\nu}\in\mathbb{C}$ for the
solution $\widehat{u}_{m,j}$ to (\ref{SRCc}) with%
\begin{equation}
\widehat{u}_{m,j}\left(  r\right)  =A_{m,j,\nu}\,f_{m,\nu}^{(1)}\left(  sr\right)
+B_{m,j,\nu}\,f_{m,\nu}^{(2)}\left(  sr\right)  \qquad\text{for all
}r>R.\label{uhatlincomb}%
\end{equation}
Standard asymptotic expansions of modified Bessel functions
\cite[\S 10.40.1--4]{NIST:DLMF} for large argument provide%
\begin{equation}%
\begin{array}
[c]{ll}%
f_{m,\nu}^{(1)}\left(  z\right)  =\frac{\operatorname*{e}^{-z}}{z^{\nu+1/2}%
}\left(  1+\mathcal{O}\left(  z^{-1}\right)  \right)  ,\quad & \frac{d}{dz}%
f_{m,\nu}^{(1)}\left(  z\right)  =-\frac{\operatorname*{e}^{-z}}{z^{\nu+1/2}%
}\left(  1+\mathcal{O}\left(  z^{-1}\right)  \right)  ,\\
f_{m,\nu}^{(2)}\left(  z\right)  =\frac{\operatorname*{e}^{z}}{z^{\nu+1/2}%
}\left(  1+\mathcal{O}\left(  z^{-1}\right)  \right)  ,\quad & \frac{d}{dz}%
f_{m,\nu}^{(2)}\left(  z\right)  =\frac{\operatorname*{e}^{z}}{z^{\nu+1/2}%
}\left(  1+\mathcal{O}\left(  z^{-1}\right)  \right)
\end{array}
\label{fundsysasy}%
\end{equation}
with the Landau-$\mathcal{O}$ notation as $\mathbb{C}_{\geq0}\ni
z\rightarrow\infty$. Consequently, for any fixed $s\in\CutC$ it holds%
\begin{equation}
r^{\nu+1/2}\left(  \frac{d}{dr}f_{m,\nu}^{\left(  j\right)  }\left(
sr\right)  +sf_{m,\nu}^{\left(  j\right)  }(sr)\right)  =%
\begin{cases}
\mathcal{O}\left(  \frac{1}{r}\right)   & \text{if }j=1,\\
\mathcal{O}\left(  \operatorname*{e}{}^{sr}\right)   & \text{if }j=2
\end{cases}
\quad\text{as }r\rightarrow\infty\label{fundsrad}%
\end{equation}
and this implies that the function $u$ given by (\ref{uhutfourier}) satisfies
condition $\operatorname*{IC}_{s}^{\operatorname*{Green}}$ if and only if
$B_{m,j,\nu}=0$ in (\ref{uhatlincomb}) for all $m\in\mathbb{N}_{0}$ and $j\in
J_{m}$. 
Moreover, (\ref{fundsrad}) and the left relation for $f_{m,\nu
}^{\left(  1\right)  }$ in (\ref{fundsysasy}) reveal the implication%
\begin{equation}
\operatorname*{IC}\nolimits_{s}^{\operatorname*{Green}}\implies
\operatorname*{IC}\nolimits_{s}^{\operatorname*{strong}}.\label{GreenStrong}%
\end{equation}
This and the Dirichlet boundary condition $\widehat{u}(R,\xi)=u\circ\psi
(R,\xi)=g\circ\psi(R,\xi)$ at $r=R$ imply
\[
A_{m,j,\nu}=\frac{1}{f_{m,\nu}^{(1)}\left(  sR\right)  }\widehat{g}_{m,j}%
\]
for the Fourier coefficients $\widehat{g}_{m,j}$ of the Dirichlet data $g\in
H^{1/2}\left(  S_{R}\right)  $ given by
\begin{equation}
\widehat{g}(\xi)=\sum_{m\in\mathbb{N}_{0}}\sum_{j\in J_{m}}\widehat{g}_{m,j}%
Y_{m,j}\left(  \xi\right)  . \label{defghat}%
\end{equation}
Hence, the solution operator $\operatorname*{S}_{\operatorname*{D}}\left(s\right)  $
for the exterior Helmholtz problem\ (\ref{eqn:exterior_problem})
maps the Dirichlet data $\widehat{g}\equiv g\in H^{1/2}\left(  S_{R}\right)  $ in
the form\ (\ref{defghat}) to%
\begin{equation}
\widehat{u}\left(  r,\xi\right)  =\sum_{m\in\mathbb{N}_{0}}\frac{f_{m,\nu}%
^{(1)}\left(  sr\right)  }{f_{m,\nu}^{(1)}\left(  sR\right)  }\sum_{j\in
J_{m}}\widehat{g}_{m,j}Y_{m,j}\left(  \xi\right)  . \label{specsolrep}%
\end{equation}
Its normal (radial) derivative results for all $s\in\CutC$
in the $\operatorname*{DtN}\left(
s\right)  $ operator
\begin{equation}
\operatorname*{DtN}\left(  s\right)  \,g:=\frac{1}{R}\sum_{m\in\mathbb{N}_{0}%
}z_{m,\nu}\left(  sR\right)  \sum_{j\in J_{m}}\widehat{g}_{m,j}Y_{m,j}\left(
\xi\right)  \quad\text{with}\quad z_{m,\nu}\left(  s\right)  :=s\frac{\frac
{d}{ds}f_{m,\nu}^{(1)}\left(  s\right)  }{f_{m,\nu}^{(1)}\left(  s\right)
}. \label{eqn:DtN_def}%
\end{equation}

\subsection{Analysis of the spectral coefficients}

\label{sub:Analysis of the spectral coefficients} The following theorem
generalises and improves the results from \cite{MelenkSauterMathComp,Nedelec01}
for purely imaginary wavenumbers $s\in\operatorname{i}\mathbb{R}$ to general
wavenumbers $s\in\mathbb{C}_{\geq0}$ and spatial dimensions $n$.

\begin{theorem}
[characterisation of $z_{m,\nu}$]\label{thm:zmn} For any $m\in\mathbb{N}_{0}$
and $\nu\in\tfrac{1}{2}\mathbb{N}_{0}$, the spectral
coefficients $z_{m,\nu}$
are holomorphic in $\CutC$ with continuous
extension to $\mathbb{C}_{\geq0}$ and

\begin{enumerate}[label=(\roman*)]
\item $\displaystyle\overline{z_{m,\nu}\left(  s\right)  }=z_{m,\nu}(\overline{s})$ for
all $s\in\mathbb{C}_{\geq0}$.
\end{enumerate}
Moreover, the following holds for all $s\in\mathbb{C}_{\geq0}$ and $\left(
m,\nu\right)  \neq(0,0)$ for some constant $c_{2}>0$:
\begin{enumerate}[label=(\roman*)]\setcounter{enumi}{1}
\item $0\leq-\operatorname{Re}z_{0,0}\left(  s\right)  \leq\tfrac{1}%
{2}+\operatorname{Re}s$\qquad and\qquad$\nu+1/2\leq-\operatorname{Re}z_{m,\nu
}\left(  s\right)  \leq m+2\nu+\operatorname{Re}s$.

\item If $\pm\operatorname{Im}s>0$, then \qquad$0<\mp$%
{$\operatorname{Im}$}$z_{0,0}\left(  s\right)  \leq c_{2}+|\operatorname{Im}%
s|$ \quad and \quad$0<\mp${$\operatorname{Im}$}$z_{m,\nu}\left(  s\right)
\leq\left\vert {\operatorname{Im}}s\right\vert $.
\end{enumerate}
\end{theorem}

The following lemma prepares the proof of \Cref{thm:zmn} below. The
modulus of the Hankel function $H_{\mu}^{(1)}=J_{\mu}+\operatorname{i}Y_{\mu}$
of order $\mu\in\mathbb{R}$ reads
\begin{equation}
M_{\mu}\left(  x\right)  =\sqrt{J_{\mu}^{2}\left(  x\right)  +Y_{\mu}%
^{2}\left(  x\right)  }\qquad\forall x>0 \label{DefMmue}%
\end{equation}
in the standard notation with the Bessel functions $J_{\mu}$ and $Y_{\mu}$ of
first and second kind%
~\cite[\S10]{NIST:DLMF}.

\begin{lemma}
[estimates for Hankel modulus functions]%
\label{lem:estimates for Hankel modulus functions} The modulus
$M_{\mu}$ satisfies
\begin{align}
M_{\mu}^{2}\left(  x\right)   &  \leq-x\frac{d}{dx}\left(  M_{\mu}^{2}\left(
x\right)  \right)  \leq2\mu M_{\mu}^{2}\left(  x\right)   & \forall x &
>0,\quad\forall\mu\geq1/2,\label{eqn:modulus_bound_mu_large}\\
2\mu\,M_{\mu}^{2}\left(  x\right)   &  \leq-x\frac{d}{dx}\left(  M_{\mu}%
^{2}\left(  x\right)  \right)  \leq M_{\mu}^{2}\left(  x\right)   & \forall x
&  >0,\quad\forall\mu\in\left[  0,\frac{1}{2}\right]
.\label{eqn:modulus_bound_mu_small}%
\end{align}
(In particular, $M_{1/2}^{2}\left(  x\right)  =-x\frac{d}{dx}\left(
M_{1/2}^{2}\left(  x\right)  \right)  $.) The bounds in
(\ref{eqn:modulus_bound_mu_large})--(\ref{eqn:modulus_bound_mu_small}) hold
with a strict inequality ($\leq$ replaced by $<$) whenever $\mu\not \in
\{0,1/2\}$ and are sharp in the sense that
\[
\lim_{x\rightarrow0}\left(  -x\frac{\frac{d}{dx}\left(  M_{\mu}^{2}\left(
x\right)  \right)  }{M_{\mu}^{2}\left(  x\right)  }\right)  =2\mu
\qquad\text{and\qquad}\lim_{x\rightarrow\infty}\left(  -x\frac{\frac{d}%
{dx}\left(  M_{\mu}^{2}\left(  x\right)  \right)  }{M_{\mu}^{2}\left(
x\right)  }\right)  =1.
\]

\end{lemma}

\proof
A holomorphic extension of $M_{\mu}^{2}$ in (\ref{DefMmue}) to $\mathbb{C}%
_{>0}$ is known from Nicholson's integral formula \cite[\S 13.73]%
{Watson} as%
\begin{equation}
M_{\mu}^{2}\left(  z\right)  :=J_{\mu}^{2}\left(  z\right)  +Y_{\mu}%
^{2}\left(  z\right)  =\frac{8}{\pi^{2}}\int_{0}^{\infty}\cosh\left(  2\mu
t\right)  K_{0}\left(  2z\sinh t\right)  dt\qquad\forall z\in\mathbb{C}%
_{>0}\label{eqn:modulus_int}%
\end{equation}
with the modified Bessel function $K_{0}$ of order $0$. Since $\cosh$ and
$\sinh$ are positive for positive arguments and increasing while $K_{0}\left(
x\right)  $ is decreasing for increasing $x>0$, it follows that $M_{\mu}%
^{2}(x)$ is decreasing in the real variable $x>0$ for fixed $\mu\geq0$. A
differentiation under the integral for positive $z\leftarrow x>0$ (both sides
in\ (\ref{eqn:modulus_int}) are analytic for $\operatorname{Re}z>0$) provides
\[
x\frac{d}{dx}\left(  M_{\mu}^{2}\left(  x\right)  \right)  =\frac{8}{\pi^{2}%
}\int_{0}^{\infty}2x\sinh\left(  t\right)  \cosh\left(  2\mu t\right)
K_{0}^{\prime}\left(  2x\sinh\left(  t\right)  \right)  dt.
\]
Employ $\frac{d}{dt}K_{0}\left(  2x\sinh\left(  t\right)  \right)
=2x\cosh(t)\,K_{0}^{\prime}(2x\sinh(t))$ from the chain rule and%
\[
{g_{\mu}}\left(  {t}\right)  :={\tanh}\left(  {t}\right)  {\coth}\left(
{2\mu\,t}\right)  {\ =\frac{\sinh(t)}{\sinh(2\mu\,t)}\,\frac{\cosh(2\mu
\,t)}{\cosh(t)}\ \qquad\forall t\in(0,\infty)}%
\]
to rewrite the last integral in the form
\begin{equation}
x\frac{d}{dx}\left(  M_{\mu}^{2}\left(  x\right)  \right)  =\frac{8}{\pi^{2}%
}\int_{0}^{\infty}g_{\mu}\left(  t\right)  \sinh\left(  2\mu t\right)
\frac{d}{dt}\left(  K_{0}\left(  2x\sinh\left(  t\right)  \right)  \right)
dt.\label{eqn:x_modulus_identity}%
\end{equation}
A straightforward calculation yields%
\[
g_{\mu}^{\prime}\left(  t\right)  =\frac{d_{\mu}\left(  t\right)  }{\cosh
^{2}\left(  t\right)  \sinh^{2}\left(  2\mu t\right)  }\quad\text{with\ }%
d_{\mu}\left(  t\right)  :=w\left(  2\mu t\right)  -2\mu
w\left(  t\right)
\quad\text{and}\quad
w\left(  t\right)  :=\sinh\left(  t\right)  \cosh\left(  t\right)
\]
for all $t\in\mathbb R$.
The derivative of $d_{\mu}$ satisfies%
\[
d_{\mu}^{\prime}\left(  t\right)  =2\mu\left(  \cosh\left(  4\mu t\right)
-\cosh\left(  2t\right)  \right)  
\qquad\forall t\in\mathbb R.
\]
For $t>0$, the functions $\sinh\left(  t\right)  $ and $\cosh\left(  t\right)
$ are positive and increasing so these properties are inherited by $w\left(
t\right)  $. For $\mu>1/2$ we have $d_{\mu}^{\prime}\left(  t\right)  >0$ and
$d_{\mu}\left(  0\right)  =0$ implies $d_{\mu}\left(  t\right)  >0$. This
allows us to conclude that $g_{\mu}^{\prime}\left(  t\right)  $ is positive
for $\mu>1/2$. An analogous reasoning verifies that $g_{\mu}^{\prime}\left(
t\right)  $ is negative for $\mu\in\left]  0,1/2\right[  $. 
Hence $g_{\mu}$ is monotonic and the limits $g_{\mu}\left(  t\right)  \rightarrow
1/\left(  2\mu\right)  $ as $t\rightarrow0$ and $g_{\mu}\left(  t\right)
\rightarrow1$ as $t\rightarrow\infty$ result in%
\begin{equation}
\min\left\{  1,1/\left(  2\mu\right)  \right\}  \leq g_{\mu}\left(  t\right)
\leq\max\left\{  1,1/\left(  2\mu\right)  \right\}  \quad\forall t,\mu
\in\left(  0,\infty\right)  \label{eqn:g_mu_bounds}%
\end{equation}
with strict inequalities for $\mu\neq1/2$. For $t>0$, the function $\sinh t>0$
is increasing while the function $K_{0}\left(  x\right)  $ is decreasing in
$\left(  0,\infty\right)  $. Thus the term $\frac{d}{dt}\left(
K_{0}\left(  2x\sinh\left(  t\right)  \right)  \right)  <0$ is negative. Hence
(\ref{eqn:x_modulus_identity})--(\ref{eqn:g_mu_bounds}) reveal
\[
\max\left\{  1,2\mu\right\}  T_{\mu}\left(  x\right)  \leq x\frac{d}%
{dx}\left(  M_{\mu}^{2}\left(  x\right)  \right)  \leq\min\left\{
1,2\mu\right\}  T_{\mu}\left(  x\right)  \quad\forall x,\mu\in\left(
0,\infty\right)
\]
(with strict inequality for $\mu\neq1//2$) for
\begin{equation}
T_{\mu}\left(  x\right)  :=\frac{8}{\pi^{2}}\int_{0}^{\infty}\frac
{\sinh\left(  2\mu t\right)  }{2\mu}\frac{d}{dt}\left(  K_{0}\left(
2x\sinh\left(  t\right)  \right)  \right)  dt.\label{eqn:T_bound}%
\end{equation}
The series and the asymptotic expansion%
~\cite[\S10.25.2, .31.2, .40.2]{NIST:DLMF} of $K_{0}$ result in
\[
\sin\left(  2\mu t\right)  K_{0}\left(  2x\sinh\left(  t\right)  \right)
\rightarrow0\quad\text{as }t\rightarrow0\text{ and as }t\rightarrow\infty.
\]
This leads to vanishing boundary terms in the integration by parts%
\begin{equation}
T_{\mu}\left(  x\right)  =-\frac{8}{\pi^{2}}\int_{0}^{\infty}\cosh\left(  2\mu
t\right)  K_{0}\left(  2x\sinh\left(  t\right)  \right)
dt\overset{\text{(\ref{eqn:modulus_int})}}{=}-M_{\mu}^{2}\left(  x\right).
\label{eqn:T_integration_parts}%
\end{equation}
The combination of
(\ref{eqn:T_bound})--(\ref{eqn:T_integration_parts}) proves
(\ref{eqn:modulus_bound_mu_large})--(\ref{eqn:modulus_bound_mu_small}) for
$x,\mu\in\left(  0,\infty\right)  $ (with strict inequality for $\mu\neq1/2$).
By continuity of $M_{\mu}^{2}=J_{\mu}^{2}+Y_{\mu}^{2}$ in $\mu$ ($J_{\mu}$ and
$Y_{\mu}$ are entire in $\mu$ \cite[\S 10.2(ii)]{NIST:DLMF}), 
(\ref{eqn:modulus_bound_mu_small}) also holds for $\mu=0$. The exactness of
the estimates in the limits $x\rightarrow0$ and $x\rightarrow\infty$ follows
from the series and the asymptotic expansions of $M_{\mu}^{2}\left(  x\right)
$, see, e.g., \cite[\S 10.2, 10.8, and 10.18]{NIST:DLMF} for further details.%
\endproof
\medskip

\textbf{Proof of \Cref{thm:zmn}.} Consider any $m\in\mathbb{N}_{0}$ and
set $\mu:=m+\nu$ with $\nu$ from (\ref{eqn:nu_def}).

\medskip\noindent\textbf{Step 1} departs with general remarks. The modified
Bessel function $K_{m+\nu}=K_{\mu}$ has no zeroes \cite[pp. 511--513]{Watson}
in $\mathbb{C}_{\geq0}$ and is holomorphic \cite[\S 10.25.ii]{NIST:DLMF} in
$\CutC$. Hence $z_{m,\nu}\left(  s\right)  $
is well defined and holomorphic for $s\neq0$. The product rule and the
recurrence relations \cite[\S 10.29.2]{NIST:DLMF} provide
\begin{equation}
z_{m,\nu}\left(  s\right)  =s\frac{K_{\mu}^{\prime}\left(  s\right)  }{K_{\mu
}\left(  s\right)  }-\nu=m-s\frac{K_{\mu+1}(s)}{K_{\mu}(s)}%
.\label{eqn:zmn_char}%
\end{equation}
The behavior of $K_{\mu}$ for small $s\in\CutC$ is known as
\[
K_{\mu}\left(  s\right)  =\left\{
\begin{array}
[c]{lc}%
\frac{(\mu-1)!}{2}\left(  \frac{2}{s}\right)  ^{\mu}+o\left(  s^{-\mu}\right)
& \text{if }\mu>0,\\
-\ln\left(  \frac{s}{2}\right)  -\gamma+o\left(  1\right)   & \text{if }\mu=0
\end{array}
\right.
\]
(with Euler--Mascheroni constant $\gamma=0.5772\ldots$), e.g., from
\cite[\S 10.30.2, .31.2]{NIST:DLMF}. Hence
\[
\lim_{\substack{s\in\CutC\\s\rightarrow
0}}s\frac{K_{\mu+1}\left(  s\right)  }{K_{\mu}\left(  s\right)  }=2\mu
\]
reveals that $z_{m,\nu}$ admits a continuous extension from $\CutC$ to $\mathbb{C}_{\geq0}$ by%
\begin{equation}
z_{m,\nu}\left(  0\right)  :=\lim_{\substack{s\in\CutC\\s\rightarrow0}}z_{m,\nu}\left(  s\right)  =m-2\mu=-\left(
m+2\nu\right)  .\label{limzmnue}%
\end{equation}
For any real parameter $\mu\in\mathbb{R}$, the modified Bessel function
$K_{\mu}(t)$ is real-valued for positive real arguments $t>0$
\cite[\S 10.25.ii]{NIST:DLMF} and satisfies $\overline{K_{\mu}\left(
z\right)  }=K_{\mu}\left(  \overline{z}\right)  $ \cite[\S 10.34.7]{NIST:DLMF}
for all $z\in\CutC$. Since complex conjugation
commutes with differentiation, this and the definition of $z_{m,\nu}(s)$ prove
\Part (i).

\medskip\noindent\textbf{Step 2} verifies \Part (ii)--(iii) for purely
imaginary wavenumbers $s=-\operatorname*{i}k$ with $k>0$. Recall the Bessel
functions $J_{\mu}$ and $Y_{\mu}$ of first and second kind \cite[\S 10.2]%
{NIST:DLMF}. The representation of $K_{\mu}(s)$ through
Hankel functions $H_{\mu}^{(1)}(k)=J_{\mu}(k)+\operatorname{i}Y_{\mu}%
(k)$~\cite[\S 10.27.8]{NIST:DLMF} is given~by%
\[
2K_{\mu}\left(  s\right)  =\pi\operatorname*{i}\nolimits^{\mu+1}{H_{\mu}%
^{(1)}}\left(  k\right)  \quad\text{with}\quad2\frac{d}{ds}K_{\mu}%
(s)=\pi\operatorname*{i}\nolimits^{\mu+2}\frac{d}{dk}{H_{\mu}^{(1)}}\left(
k\right)  \qquad\text{for }s=-\operatorname{i}k.
\]
This and (\ref{eqn:zmn_char}) result in
\begin{equation}
z_{m,\nu}\left(  -\operatorname*{i}k\right)  =k\frac{\frac{d}{dk}H_{\mu}%
^{(1)}(k)}{H_{\mu}^{(1)}(k)}-\nu=w_{\mu}(k)-\nu\qquad\text{with}\quad w_{\mu
}(k):=k\frac{\frac{d}{dk}H_{\mu}^{(1)}(k)}{H_{\mu}^{(1)}(k)}%
.\label{eqn:zmn_char_imag}%
\end{equation}
The investigation of $w_{\mu}(k)$ departs with the decomposition into its real
and imaginary parts known from~\cite[p.~1877]{MelenkSauterMathComp} and is
repeated here for convenience. First observe
\[
w_{\mu}(k)=k\frac{H_{\mu}^{(2)}(k)\frac{d}{dk}H_{\mu}^{(1)}(k)}{H_{\mu}%
^{(2)}(k)H_{\mu}^{(1)}(k)}=k\frac{J_{\mu}(k)J_{\mu}^{\prime}(k)+Y_{\mu
}(k)Y_{\mu}^{\prime}(k)+\operatorname{i}\left(  J_{\mu}(k)Y_{\mu}^{\prime
}(k)-Y_{\mu}(k)J_{\mu}^{\prime}(k)\right)  }{J_{\mu}^{2}(k)+Y_{\mu}^{2}(k)}%
\]
by expanding (\ref{eqn:zmn_char_imag}) with $H_{\mu}^{(2)}(k)=J_{\mu
}(k)-\operatorname{i}Y_{\mu}(k)$. The Wronskian identity $J_{\mu}(k)Y_{\mu
}^{\prime}(k)-Y_{\mu}(k)J_{\mu}^{\prime}(k)=2/(\pi k)$ from~\cite[\S 10.5.2]%
{NIST:DLMF} is a consequence of the recurrence relations for Bessel functions.
This reveals the characterisation \cite[Eqn.~(3.10)]{MelenkSauterMathComp} of
$w_{\mu}$, namely
\[
w_{\mu}\left(  k\right)  =\frac{k}{2}\frac{\frac{d}{dk}\left(  M_{\mu}%
^{2}\left(  k\right)  \right)  }{M_{\mu}^{2}\left(  k\right)  }%
+\operatorname*{i}\frac{2}{\pi M_{\mu}^{2}\left(  k\right)  }\qquad\text{for
}M_{\mu}\left(  k\right)  =\left\vert H_{\mu}^{\left(  1\right)  }\left(
k\right)  \right\vert \text{ as in (\ref{DefMmue}).}%
\]
Hence the real and imaginary parts of $z_{m,n}(-\operatorname{i}k)$ satisfy
by~(\ref{eqn:zmn_char_imag}) that
\begin{align}
-\operatorname{Re}z_{m,\nu}\left(  -\operatorname{i}k\right)   &
=\nu-\operatorname{Re}w_{\mu}\left(  k\right)  =\nu-\frac{k}{2}\frac{\frac
{d}{dk}\left(  M_{\mu}^{2}\left(  k\right)  \right)  }{M_{\mu}^{2}\left(
k\right)  },\label{eqn:Re_z}\\
\operatorname{Im}z_{m,\nu}\left(  -\operatorname{i}k\right)   &
=\operatorname{Im}w_{\mu}\left(  k\right)  =\frac{2}{\pi M_{\mu}^{2}%
(k)}.\label{eqn:Im_z}%
\end{align}
Recall the assumption $k>0$ in this step. Lemma
\ref{lem:estimates for Hankel modulus functions} controls the real part of
$w_{\mu}(k)$ by
\begin{align}
\frac{1}{2} &  \leq-\operatorname{Re}w_{\mu}(k)\leq\mu\qquad\text{for all }%
\mu\geq1/2,\label{eqn:Re_w}\\
0 &  \leq-\operatorname{Re}w_{0}(k)\leq1/2.\label{eqn:Re_w_0}%
\end{align}
Since $kM_{\mu}^{2}(k)$ is non-increasing for $\mu\geq1/2$ and converges to
$2/\pi$ as $k\rightarrow\infty$ \cite[\S 13.74]{Watson}, one has $2\leq\pi
kM_{\mu}^{2}(k)<\infty$ and
\[
0<\operatorname{Im}w_{\mu}(k)=\frac{2}{\pi M_{\mu}^{2}(k)}\leq k\qquad
\text{for all }\mu\geq1/2.
\]
The situation is different for $\mu=0$ where $M_{0}^{2}(k)$ admits no lower
bound for $k\rightarrow0$ that is linear in $1/k$. However, the asymptotic
expansion~\cite[\S 13.75]{Watson} of $M_{0}^{2}(k)$ reveals the convergence
$kM_{0}^{2}(k)\rightarrow2/\pi$ as $k\rightarrow\infty$ with the lower bound
(see \cite[\S 10.18.iii]{NIST:DLMF})%
\begin{equation}
\frac{2}{\pi}\left(  1-\frac{1}{8k^{2}}\right)  \leq kM_{0}^{2}\left(
k\right)  \qquad\text{for all }k>0.\label{eqn:Im_w}%
\end{equation}
Fix $k_{0}>2^{-3/2}$ and set $C_{0}:=\frac{k_{0}}{8k_{0}^{2}-1}$. This leads
to%
\begin{equation}
\frac{2}{\pi M_{0}^{2}(k)}\leq k\frac{8k^{2}}{8k^{2}-1}=k+\frac{k}{8k^{2}%
-1}\leq k+C_{0}\qquad\forall k>k_{0}.\label{eqn:Im_z_0_tmp}%
\end{equation}
On the other hand, we have already seen that $M_{0}^{2}$ is non-increasing so
that%
\[
\sup_{0<k<k_{0}}\frac{2}{\pi M_{0}^{2}(k)}\leq\frac{2}{\pi M_{0}^{2}(k_{0}%
)}=:C_{0}^{\prime}.
\]
For $c_{2}:=C_{0}+C_{0}^{\prime}$, this establishes
\begin{equation}
0<\operatorname{Im}w_{0}(k)=\frac{2}{\pi M_{0}^{2}(k)}\leq c_{2}+k\qquad\text{
}\forall k>0.\label{eqn:Im_w0}%
\end{equation}
The combination of\ (\ref{eqn:Re_z})--(\ref{eqn:Im_z}) with (\ref{eqn:Re_w})
and (\ref{eqn:Im_w})\ reads
\begin{equation}
\nu+\frac{1}{2}\leq-\operatorname{Re}\,z_{m,\nu}\left(  -\operatorname{i}%
k\right)  \leq\nu+\mu\qquad\text{and}\qquad0<{\operatorname{Im}}\,z_{m,\nu
}\left(  -\operatorname{i}k\right)  \leq k\label{eqn:z_bound}%
\end{equation}
for all $\mu\geq1/2$, whereas (\ref{eqn:Re_w_0}) and (\ref{eqn:Im_w0}) result
for $\mu=0=m=\nu$ in
\begin{equation}
0\leq-\operatorname{Re}\,z_{0,0}\left(  -\operatorname{i}k\right)  \leq
\frac{1}{2}\qquad\text{and}\qquad0<{\operatorname{Im}}\,z_{0,0}\left(
-\operatorname{i}k\right)  \leq c_{2}+k.\label{eqn:z_bound_0}%
\end{equation}

\medskip\noindent\textbf{Step 3} extends the analysis to general wavenumbers
$s\in{\mathbb{C}}_{\geq0}$.
The behaviour of the modified Bessel functions for large $|s|$ and any fixed $\mu$ 
is known \cite[\S 10.40.10]{NIST:DLMF} as
\[
K_{\mu}(s)=\sqrt{\frac{\pi}{2s}}\operatorname*{e}\nolimits^{-s}\left(
1+\frac{a_{1}(\mu)}{s}+R_{2}(\mu,s)\right)
\]
where $a_{1}(\mu)=(4\mu^{2}-1)/8$ and the remainder $R_{2}(\mu,s)\leq
C|s|^{-2}$ uniformly tends to zero for large $s\in\mathbb{C}_{\geq0}$.
Hence~(\ref{eqn:zmn_char}) reveals, for any $\mu=m+\nu\in
\tfrac{1}{2}\mathbb{N}_{0}$, the asymptotic behaviour
\[
\lim_{\substack{s\in\mathbb{C}_{\geq0}\\s\rightarrow\infty}}\left(
m-z_{m,\nu}\left(  s\right)  -s\right)  =\lim_{\substack{s\in\mathbb{C}%
_{\geq0}\\s\rightarrow\infty}}s\frac{K_{\mu+1}\left(  s\right)  -K_{\mu}%
(s)}{K_{\mu}\left(  s\right)  }=a_{1}(\mu+1)-a_{1}(\mu)=\mu+\frac{1}{2}.
\]
In other words, for any $\varepsilon>0$ there is some $r>0$ such that all
$s\in\mathbb{C}_{\geq0}$ with $|s|=r$ satisfy
\begin{align*}
\nu+\frac{1}{2}-\varepsilon &  \leq-\operatorname{Re}\,(z_{m,\nu}%
(s)+s)\leq-\operatorname{Re}\,z_{m,\nu}(s)\leq\nu+\frac{1}{2}%
+\operatorname{Re}\,s+\varepsilon,\\
-\varepsilon &  \leq\mp\operatorname{Im}\,(z_{m,\nu}(s)+s)\leq\mp
\operatorname{Im}\,z_{m,\nu}(s)\leq|\operatorname{Im}\,s|+\varepsilon
\qquad\text{whenever }0\leq\pm \operatorname{Im} s.
\end{align*}
Since Step 2 controls the real part of
$z_{m,\nu}$ on the imaginary axis (by (\ref{eqn:z_bound}%
)--(\ref{eqn:z_bound_0}) and $z_{m,\nu}(\overline{s})=\overline{z_{m,\nu
}\left(  s\right)  }$ from \Part (i)), this leads for all $s\in\partial
(B_{r}\cap\mathbb{C}_{\geq0})$ to%
\begin{align}
\nu+\frac{1}{2}-\varepsilon &  \leq-\operatorname{Re}\,z_{m,\nu}(s)\leq\nu
+\mu+\operatorname{Re}\,s+\varepsilon\qquad\text{if }1/2\leq\mu=m+\nu
,\label{eqn:Rz_eps}\\
0 &  \leq-\operatorname{Re}\,z_{0,0}(s)\leq\frac{1}{2}+\operatorname{Re}%
\,s+\varepsilon.\label{eqn:Rz_eps_0}%
\end{align}
For the imaginary part, where $s\in\partial Q_{r}^{+}$ varies on the boundary
of $Q_{r}^{+}:=\{s\in\overline{B_{r}}\cap\mathbb{C}_{\geq0}%
\ :\ \operatorname{Im}s\geq0\}$ in the first quadrant,
it follows from (\ref{eqn:z_bound})--(\ref{eqn:z_bound_0}) with \Part (i) that
\begin{equation}
-\varepsilon\leq-\operatorname{Im}z_{m,\nu}(s)\leq c_{n}+\operatorname{Im}%
s+\varepsilon\label{eqn:Iz_eps}%
\end{equation}
with $c_{2}$ from~(\ref{eqn:Im_w0}) and $c_{n}:=0$ for all $n\geq3$.
It is well-known that the real and imaginary parts of the holomorphic function
$z_{m,\nu}(s)$ (by Step 1) are harmonic. The maximum principle for harmonic
functions \cite[Thm.~2.2]{Gilbarg} reveals\ (\ref{eqn:Rz_eps}%
)--(\ref{eqn:Rz_eps_0}) for all $s\in B_{r}\cap\mathbb{C}_{\geq0}$ and
(\ref{eqn:Iz_eps}) for all $s\in Q_{r}^{+}$. Since this holds for all
$\varepsilon>0$ and sufficiently large $r>0$, this concludes the proof of \Part
(ii) for the real part and leads for all $r>0$ to
\[
-\varepsilon\leq-\operatorname{Im}z_{m,\nu}(s)\leq c_{n}+\operatorname{Im}%
s+\varepsilon\qquad\text{for all }s\in Q_{r}^{+}.
\]
The strict sign $0<\mp\operatorname{Im}z_{m,\nu}(s)$ on the imaginary axis
with $s\in\operatorname{i}\mathbb{R}\backslash\left\{  0\right\}  $
(by~(\ref{eqn:z_bound})--(\ref{eqn:z_bound_0}) and \Part (i))
contradicts that $\operatorname{Im}z_{m,\nu}(s)\equiv0$ vanishes in $Q_{r}%
^{+}$. Hence $0<-\operatorname{Im}z_{m,\nu}(s)$ also holds for all
$s\in\mathrm{int}\,Q_{r}^{+}$ in the interior of $Q_{r}^{+}$ by the (strong)
maximum principle. This and $\operatorname{Im}z_{m,\nu}(\overline
{s})=-\operatorname{Im}z_{m,\nu}(s)$ from \Part (i) conclude the proof of \Part
(iii).\qed%
\endproof

\subsection{Spectral representation of the Laplace DtN}\label{sub:DtN_Laplace}

The discussion in Section
\ref{sub:Spectral expansion of the Helmholtz DtN operator} reveals for $s=0$
that any $u\in H_{\operatorname*{loc}}^{1}(B_{R}^{+})$ with a series
expansion of the form (\ref{uhutfourier}) solves the Laplace equation $-\Delta
u=0$ if its Fourier coefficients $\widehat{u}_{m,j}$ given by
(\ref{uhutfourier}) solve the ODE
\begin{equation}
-\frac{1}{r^{n-1}}\partial_{r}\left(  r^{n-1}\partial_{r}\widehat{u}_{m,j}\left(
r\right)  \right)  +\frac{m\left(  m+n-2\right)  }{r^{2}}\widehat{u}_{m,j}\left(
r\right)  =0\label{locPDEdelta}%
\end{equation}
for all $r>R$ and $m\in\mathbb{N}_{0},j\in J_{m}$. A fundamental system for this ODE is
given by%
\[
f_{m,\nu}^{(1)}\left(  r\right)  =r^{-m-2\nu}\qquad\text{and}\qquad f_{m,\nu
}^{(2)}\left(  r\right)  :=\left\{
\begin{array}
[c]{ll}%
r^{m} & \text{if }m+2\nu>0,\\
\log r & \text{if }m+2\nu=0
\end{array}
\right.  \qquad\text{for all }m\in\mathbb{N}_{0}.
\]
From this, we deduce that the behaviour of $f_{m,\nu}^{(2)}$ contradicts the
condition $\operatorname*{IC}_{0}^{\operatorname*{strong}}$ while $f_{m,\nu
}^{\left(  1\right)  }$ satisfies $\operatorname*{IC}_{0}%
^{\operatorname*{strong}}$. Hence, all coefficients $B_{m,\nu,j}$ in
(\ref{uhatlincomb}) must be zero. By matching the other coefficients with the
series expansion of the Dirichlet data $g$ (with~\eqref{defghat}), we obtain%
\begin{equation}
\widehat{u}\left(  r,\xi\right)  =\sum_{m\in\mathbb{N}_{0}}\frac{f_{m,\nu
}^{(1)}\left(  r\right)  }{f_{m,\nu}^{(1)}\left(  R\right)  }\sum_{j\in J_{m}%
}\widehat{g}_{m,j}Y_{m,j}\left(  \xi\right)  =\sum_{m\in\mathbb{N}_{0}}\left(
\frac{R}{r}\right)  ^{m+2\nu}\sum_{j\in J_{m}}\widehat{g}_{m,j}Y_{m,j}\left(
\xi\right)  \quad\forall r>R,\xi\in S_{1}.\label{uhatexpLap}%
\end{equation}
Its normal (radial) derivative results in the $\operatorname*{DtN}\left(
s\right)  $ operator for the Laplacian%
\[
{\DtNLaplace}g=\frac{1}{R}\sum_{m\in\mathbb{N}_{0}%
}z_{m,\nu}^{\Delta}\sum_{j\in J_{m}}\widehat{g}_{m,j}Y_{m,j}\left(
\xi\right)  \quad\text{with\quad}z_{m,\nu}\left(0\right):=R\frac{\frac{d}%
{dR}f_{m,\nu}^{(1)}\left(  R\right)  }{f_{m,\nu}^{(1)}\left(  R\right)
}=-(m+2\nu).
\]
It follows from (\ref{limzmnue}) in Step 1 of the proof of Theorem
\ref{thm:zmn} that $z_{m,\nu}(0)$ is the limit $z_{m,\nu}(s)$ as
$s\rightarrow0$. \W{This justifies the strong convergence $\operatorname{DtN}(s)\to \operatorname{DtN}(0)$ in \Cref{thm:DtN}.}
Next we prove \Cref{lem:two dimensions}.

\medskip
\begin{proof}[Proof of \Cref{lem:two dimensions}]
Consider\ problem
(\ref{eqn:exterior_problem}) with $\operatorname*{IC}_{0}%
^{\operatorname*{Green}}$ for the Laplacian ($s=0$) and $n=2$, $R\neq1$. 
Since $f_{0,0}^{(2)}$ (and not
$\W{f_{0,0}^{(1)}\equiv 1}$) satisfies $\operatorname*{IC}_{0}^{\operatorname*{Green}}$,
only the first term $m=0$ (with $J_{0}=\left\{  1\right\}  $) in the
series representation (\ref{uhatexpLap}) of the solution changes by 
\W{the substitution of the constant
\begin{align*}
	1=\frac{f_{0,0}^{(1)}\left(  r\right)  }{f_{0,0}^{(1)}\left(  R\right)  }
	\to 
	\frac{f_{0,0}^{(2)}\left(  r\right)  }{f_{0,0}^{(2)}\left(
	R\right)  }=\frac{\log r}{\log R}
\end{align*}
with a logarithmic term in $r$}.
This leads to the representation
\begin{equation}
\widehat{u}_{\mathrm{\log}}\left(  r,\xi\right)  =\frac{\log\left(  r\right)
}{\log\left(  R\right)  }\widehat{g}_{0,1}Y_{0,1}\left(  \xi\right)
+\sum_{m\in\mathbb{N}}\frac{f_{m,\nu}^{(1)}\left(  r\right)  }{f_{m,\nu}%
^{(1)}\left(  R\right)  }\sum_{j\in J_{m}}\widehat{g}_{m,j}Y_{m,j}\left(
\xi\right)  \label{uhatlog}%
\end{equation}
of the unique solution of Laplace problem (\ref{eqn:exterior_problem}) for
$s=0$ with $\operatorname*{IC}_{0}^{\operatorname*{Green}}$. Since
$Y_{0,1}(\xi)=Y_{0,1}$ is constant, the coefficient is explicitly given by the
integral mean of $g$, namely
\begin{equation}
\widehat{g}_{0,1}Y_{0,1}=\frac{Y_{0,1}^{2}\int_{S_{R}}g}{\Vert Y_{0,1}%
\Vert_{L^{2}\left(  S_{R}\right)  }^{2}}=\frac{1}{|S_{R}|}\int_{S_{R}%
}g.\label{meanvalueg0}%
\end{equation}
Hence \Cref{PropProblems}.2b reveals for any
$R\neq1$ that the function (\ref{uhatlog}) is the unique solution
${\Slog}  g\ $of the exterior
Laplace problem with $\operatorname*{IC}_{0}^{\operatorname*{Green}}$.
In the case $R=1$, the first summand in (\ref{uhatlog}) is only well
defined (as zero) provided $\int_{S_{1}}g=0$.
Conversely, if $\int_{S_{1}}g\ne0$ for $n=2$ and $R=1$, then~\eqref{eqn:exterior_problem} does \emph{not}
admit a solution with $\operatorname*{IC}_{0}^{\operatorname*{Green}}$.

It remains to consider the difference ${\SLaplace}g-{\Slog}  g$ which
can be written as%
\[
{\SLaplace}g-{\Slog}  g=\widehat{g}_{0,1}%
Y_{0,1}-\frac{\log\left(  r\right)  }{\log\left(  R\right)  }\widehat{g}%
_{0,1}Y_{0,1}=-\frac{\log\frac{r}{R}}{\log R}\frac{1}{|S_{R}|}\int_{S_{R}}g.
\]
Hence, the difference of the $\operatorname*{DtN}$ operators satisfies%
\[
{\DtNLaplace}g-\DtNlog
g=-\frac{1}{R\log R}\frac{1}{|S_{R}|}\int_{S_{R}}g.
\]
This finishes the proof of \Cref{lem:two dimensions}.
\end{proof}
\subsection{Proofs of Theorems 3.4 and 3.5}\label{sub:Proof of DtN}
The sharp bounds in \Cref{thm:zmn} and the explicit representation of the 
DtN operator~\eqref{eqn:DtN_def}
enable the proofs of \Cref{thm:DtN,thm:C_Friedrich} stated in \Cref{sub:DtN_Helmholtz}.

\begin{proof}[Proof of \Cref{thm:DtN}]
This proof considers the unit disc $S_{1}$ ($R=1$) and the general case $R>0$
follows from the scaling of the norms and the arising constants then also depend
on $R$. Any $g,h\in H^{1/2}(S_{1})$ with Fourier coefficients $\widehat
{g}_{m,j},\widehat{h}_{m,j}$ from~(\ref{defghat}) with respect to the $L^{2}%
(S_{1})$ orthonormal basis $\{Y_{m,j}\}_{m\in\mathbb{N}_{0},j\in J_{m}}$
result in
\begin{equation}
\left\langle \operatorname*{DtN}\left(  s\right)  g,\overline{h}\right\rangle
_{S_{1}}=\sum_{m\in\mathbb{N}_{0}}\sum_{j\in J_{m}}z_{m,\nu}(s)\,\widehat{g}%
_{m,j}\,\overline{\widehat{h}_{m,j}}\label{eqn:DtN_gh}%
=
\left\langle g,\operatorname*{DtN}\left(  s\right)  \overline{h}\right\rangle.
\end{equation}
Since $z_{m,\nu}(s)=\overline{z_{m,\nu}(\overline{s})}$ by Theorem
\ref{thm:zmn}.i, this proves \Cref{thm:DtN}.i.
Recall from the definition of the Fourier coefficients~(\ref{defghat}) and the
equivalence of norms \cite[Rem.~7.6]{LionsMagenesI} in $H^{1/2}(S_{1})$ that
there exists a constant $C_1(n)>0$ that exclusively depends on $n$ such that
\begin{equation}
\Vert g\Vert_{L^{2}(S_{1})}^{2}=\sum_{m\in\mathbb{N}_{0}}\sum_{j\in J_{m}%
}\left\vert \widehat{g}_{m,j}\right\vert ^{2}\quad\text{and}\quad\sum
_{m\in\mathbb{N}_{0}}\sum_{j\in J_{m}}\left(  1+m\right)  \left\vert \widehat
{g}_{m,j}\right\vert ^{2}\leq C_1\left(  n\right)  \left\Vert g\right\Vert
_{H^{1/2}\left(  S_{1}\right)  }^{2}.\label{eqn:Fourier_norm}%
\end{equation}
\Cref{thm:zmn}.ii--iii implies $\left\vert \operatorname{Re}z_{m,\nu
}\left(  s\right)  \right\vert \leq A+\operatorname{Re}s$ and $\left\vert
\operatorname{Im}z_{m,\nu}\left(  s\right)  \right\vert \leq B+\left\vert
\operatorname{Im}s\right\vert $ for $A=\max\left\{  \frac{1}{2},m+2\nu
\right\}  >0$ and $B=c_{2}>0$ such that
\[
\left\vert z_{m,\nu}\left(  s\right)  \right\vert ^{2}\leq\left(
A+\operatorname{Re}s\right)  ^{2}+\left(  B+|\operatorname{Im}s|\right)
^{2}=\left\vert A+\operatorname{Re}s+\operatorname*{i}\left(
B+|\operatorname*{Im}s|\right)  \right\vert ^{2}=\left\vert A\pm \operatorname*{i}%
B+s\right\vert ^{2}.%
\]
A triangle inequality results in%
\begin{align}\nonumber
\left\vert z_{m,\nu}\left(  s\right)  \right\vert  &  \leq\left\vert
A+\operatorname*{i}B+s\right\vert \leq\left\vert A\right\vert +\left\vert
B\right\vert +\left\vert s\right\vert \\
&
\leq m+2\nu+\frac{1}{2}+c_{2}+\left\vert
s\right\vert 
\label{eqn:z_m_n_absolute_bound}
  \leq C_{2}\left(  n\right)  \left(  m+1\right)  +\left\vert s\right\vert
\end{align}
for all $s\in\mathbb{C}_{\geq0}$
with $C_{2}\left(  n\right)  =1+2\nu+c_{2}>0$ that only depends on $n=2\nu+2$.

Hence (\ref{eqn:DtN_gh}) and a Cauchy inequality verify for all $g,h\in
H^{1/2}(S_{1})$ that
\begin{align*}
|\left\langle \operatorname*{DtN}\left(  s\right)  g,\overline{h}\right\rangle
_{S_{1}}|  &  \leq\sum_{m\in\mathbb{N}_{0}}\sum_{j\in J_{m}}\left(
C_{2}\left(  n\right)  \left(  1+m\right)  +\left\vert s\right\vert \right)
\,\left\vert \widehat{g}_{m,j}\right\vert \,\left\vert \widehat{h}_{m,j}\right\vert \\
&  \leq|s|\Vert g\Vert_{L^{2}(S_{1})}\left\Vert h\right\Vert _{L^{2}(S_{1}%
)}+C_1(n)C_{2}\left(  n\right)  \Vert g\Vert_{H^{1/2}(S_{1})}\left\Vert
h\right\Vert _{H^{1/2}(S_{1})}.
\end{align*}
This and a Cauchy inequality lead to the boundedness of $\operatorname*{DtN}%
\left(  s\right)  \in L(H^{1/2}(S_{1});H^{-1/2}(S_{1}))$ in the form of
\Cref{thm:DtN}.\W{ii} with $C(n)=\max\{1,C_1(n)C_{2}\left(  n\right)  \}$.
The \W{strong continuity of $\operatorname{DtN}(s)$} in the wavenumber $s\in\mathbb{C}_{\geq0}$ follows from Lebesgue's dominated convergence theorem.
\W{Indeed,~\eqref{eqn:DtN_gh}--\eqref{eqn:Fourier_norm} imply for any $s,s_0\in\CutC$ and the dual norm in $H^{-1/2}(S_1)=(H^{1/2}(S_1))^*$ that
	\begin{align}\label{eqn:DtN_converge}
		\left\Vert (\operatorname*{DtN}(s)-\operatorname*{DtN}(s_0))g\right\Vert_{H^{-1/2}(S_1)}^2
		\leq C_{1}(n)\sum
_{m\in\mathbb{N}_{0}}\sum_{j\in J_{m}}\frac{|z_{m,\nu}(s)-z_{m,\nu}(s_0)|^2}{1+m}  \left\vert \widehat
{g}_{m,j}\right\vert ^{2}.
	\end{align}
	As $s\to s_0$, the spectral coefficients $|z_{m,\nu}(s)-z_{m,\nu}(s_0)|^2/(1+m)\to 0$ remain bounded with~\eqref{eqn:z_m_n_absolute_bound} up to a constant by $(1+m)$.
	Hence using~\eqref{eqn:Fourier_norm}, the right-hand side of~\eqref{eqn:DtN_converge} is absolutely convergent and $\operatorname{DtN}(s)g\to \operatorname{DtN}(s_0)g$ in $H^{-1/2}(S_1)$ follows.} 
	\Cref{thm:zmn} controls the real and imaginary parts
of~(\ref{eqn:DtN_gh}) for $g=h$ and results with (\ref{eqn:Fourier_norm}) \W{and $c\left(  n\right)
=\tfrac{n-2}{2}$} in
\begin{align*}
c(n)\Vert g\Vert_{L^{2}(S_{1})}^{2}  &  \leq-\operatorname{Re}\left(
\left\langle \operatorname*{DtN}\left(  s\right)  g,\overline{g}\right\rangle
_{S_{1}}\right)  \leq\operatorname{Re}s\Vert g\Vert_{L^{2}(S_{1})}%
^{2}+C(n)\Vert g\Vert_{H^{1/2}(S_{1})}^{2},\\
0  &  <\mp\operatorname{Im}\left(  \left\langle \operatorname*{DtN}\left(
s\right)  g,\overline{g}\right\rangle _{S_{1}}\right)  \leq
C(n)|\operatorname{Im}s|\Vert g\Vert_{L^{2}(S_{1})}^{2}\qquad\text{if }g\neq0
\end{align*}
for some constant $C(n)$ that exclusively depends on $n$. This is \Part
(ii)--(iii) in \Cref{thm:DtN} and concludes the proof.
\end{proof}

\begin{proof}[Proof of \Cref{thm:C_Friedrich}]
	Recall the Hilbert space $H^1(\Omega,S_R)$ from~\W{\Cref{ChapPrelim}} (for $\Gamma=S_R$).

Assume for a contradiction to (\ref{eqn:C_F}) that there exists two sequences
$(v_{k})_{k\in\mathbb{N}}\subset H^1(\Omega,S_R)$ and 
$(s_k)_{k\in\mathbb N}\subset \mathbb{C}_{\geq0}$ 
with $\partial_r v_k =\operatorname{DtN}(s_k)v_k$  on $S_R$ such that
$\left\Vert v_{k}\right\Vert
_{L^{2}\left(  \Omega\right)  }=1$ for all $k\in\mathbb{N}$ and $\Vert\nabla
v_{k}\Vert_{L^{2}\left(  \Omega\right)  }+\left\Vert \partial_{r}%
v_{k}\right\Vert _{H^{-1/2}(S_R)}\rightarrow0$ as $k\rightarrow
\infty$. Since the sequence $(v_k)_{k\in\mathbb{N}}$ is bounded in $H^1(\Omega,S_R)$, 
it admits a weakly convergent
subsequence with weak limit $u\in H^1(\Omega,S_R)$. 
A standard compactness argument with the
compact embedding $H^{1}\left(  \Omega\right)  \hookrightarrow L^{2}\left(
\Omega\right)  $ and the weak lower semicontinuity of norms
reveal $\left\Vert u\right\Vert
_{L^{2}\left(  \Omega\right)  }=1$ and $\left\Vert \nabla u\right\Vert
_{L^{2}\left(  \Omega\right)  }=0
=\left\Vert \partial_{r}u\right\Vert _{H^{-1/2}\left(  S_{R}\right)}$. 
Hence $u\in H^1(\Omega,S_R)$ is constant on the 
(in particular connected) domain $\Omega$.

The traces of $v_k$ and the constant $u$ on $S_R$ have the Fourier expansions
\begin{align*}
	v_k|_{S_R}= \sum^{}_{m\in\mathbb{N}_0} \sum_{j\in J_m} \widehat{v}_{m,j}^{(k)}Y_{m,j}
	\qquad\text{and}\qquad
	u|_{S_R} = \widehat u_{0,0} Y_{0,0}
	\quad\text{for }\widehat u_{0,0}=\frac{u}{Y_{0,0}} = \frac{u}{|S_R|^{1/2}}
\end{align*}
with $\widehat{v}_{\bullet,\bullet}^{(k)}\in\mathbb{C}$
for all $k\in\mathbb{N}$ in terms of 
the $L^2(S_R)$-orthonormal Laplace-Beltrami eigenfunctions 
$Y_{\bullet,\bullet}$ from~\eqref{defiotam}.
The DtN operator representation~\eqref{eqn:DtN_def} and
$\partial_r v_k =\operatorname{DtN}(s_k) v_k$
result~in
\begin{align*}
	\frac{z_{0,\nu}(s_kR)}{R}\widehat{v}_{0,0}^{(k)}\,\overline{\widehat u_{0,0}}
	=\left\langle \operatorname{DtN}(s_k)v_k, \overline u\right\rangle_{S_R}
	=\left\langle \partial_rv_k, \overline u\right\rangle_{S_R}
	\to 0
\end{align*}
with the weak convergence $\partial_r v_k\rightharpoonup \partial_r u=0$ 
in $H^{-1/2}(S_R)$ as $k\to \infty$ in the last step.
This and $1\leq \nu+1/2\leq |z_{0,\nu}(s_k R)|$ 
for $n\geq 3$ from \Cref{thm:zmn}.ii establish 
\begin{align*}
	|\widehat{v}_{0,0}^{(k)}||\widehat u_{0,0}|
	\leq |z_{0,\nu}(s_k R)||\widehat{v}_{0,0}^{(k)}||\widehat u_{0,0}|
	= R\, 
	\left|\left\langle \operatorname{DtN}(s_k)v_k, \overline u\right\rangle_{S_R}\right|
	\to 0
	\qquad\text{as }k\to \infty.
\end{align*}
Thus either $\widehat u_{0,0}=0$ or 
$\widehat{v}_{0,0}^{(k)}\to 0$ as $k\to \infty$ hold.
The weak convergence $v_k|_{S_R}\rightharpoonup u|_{S_R}$ in $H^{1/2}(S_R)$ 
(by the continuity of the trace map) implies the convergence
$\widehat{v}_{0,0}^{(k)}\to \widehat u_{0,0}$ in $\mathbb{C}$ as $k\to \infty$
and reveals $\widehat u_{0,0}=0$ in the second case.
Hence $u=0$ follows in any case.
This contradicts $\|u\|_{L^2(\Omega)}=1$ and
concludes the proof of~\eqref{eqn:C_F}.
\end{proof}

\appendix
\section{Appendix: Proof of Proposition 3.2}\label{app:PropProblems}

The proof of \Cref{PropProblems} follows mainly from the theory in
\cite{Mclean00}, where integral equation techniques are used. The key role is
played by Green's representation formula and we sketch the arguments in the
following. The starting point is the exterior problem \cite[Eqn.\ (7.31)]{Mclean00}
( for $\mathcal{P}=-\Delta+s^{2}$, $f=0$, and $\Gamma=S_{R}$ therein)
written as:
\[
\text{find }u\in H_{\operatorname*{loc}}^{1}\left(  B_{R}^{+}\right)
\quad\text{such that\quad}\left\{
\begin{array}
[c]{rl}%
-\Delta u+s^{2}u=0 & \text{in }B_{R}^{+},\\
u=g & \text{on }S_{R},\\
\mathcal{M}u=0 & .
\end{array}
\right.
\]
Here $\mathcal{M}u=0$ is a condition for the behaviour at infinity which is
defined in \cite[(7.29)]{Mclean00}.
Theorems 8.9 and 9.6 in \cite{Mclean00} translate this condition for our setting to
\[
u\text{ satisfies }\operatorname*{IC}\nolimits_{s}^{\operatorname*{Green}}%
\iff\mathcal{M}u=0.
\]
If the exterior problem 
(\ref{eqn:exterior_problem}) with condition $\operatorname*{IC}\nolimits_{s}%
^{\operatorname*{Green}}$ has a solution $u$, then this solution satisfies
Green's representation theorem by \cite[Thm.~7.15]{Mclean00}, namely
\[
u=\operatorname*{DL}\left(  s\right)  g-\operatorname*{SL}\left(  s\right)
\psi\qquad\text{in }B_{R}^{+}%
\]
where $\operatorname*{SL}\left(  s\right)  $ and $\operatorname*{DL}\left(
s\right)  $ are the single and double layer potential operators for the
operator $\left(  -\Delta+s^{2}\right)  $ 
(see, e.g., \cite[Eqn.~(7.1)]{Mclean00} for details)
and $\psi$ is a solution of%
\begin{equation}
S\psi=\left(  -\frac{1}{2}+T\right)  g\qquad\text{on }S_{R},\label{Seq}%
\end{equation}
with the abbreviations 
$S:=\frac{1}{2}\left(  \gamma^{+}\operatorname*{SL}\left(  s\right)
+\gamma^{-}\operatorname*{SL}\left(  s\right)  \right)  $ and $T:=\frac{1}%
{2}\left(  \gamma^{+}\operatorname*{DL}\left(  s\right)  +\gamma
^{-}\operatorname*{DL}\left(  s\right)  \right)  $. 
Note that for any $g\in
H^{1/2}\left(  S_{R}\right)  $ and $h\in H^{-1/2}\left(  S_{R}\right)  $, 
$\operatorname*{DL}\left(  s\right)  g$ and $\operatorname*{SL}%
\left(  s\right)  h$ define functions in $\mathbb{R}^{n}\backslash S_{R}$ and
$\gamma^{+}$ denotes the trace from the outer domain $B_{R}^{+}$, while
$\gamma^{-}$ is the trace from the inner domain $B_{R}$. It is proved in
\cite[Theorem 6.11]{Mclean00} that the mappings
\begin{align*}
	\operatorname*{DL}\left(
		s\right)  &:H^{1/2}\left(  \Gamma\right)  \rightarrow H_{\operatorname*{loc}%
		}^{1}\left(  \mathbb{R}^{n}\backslash S_{R}\right),\\
	\operatorname*{SL}%
	\left(  s\right)  &:H^{-1/2}\left(  \Gamma\right)  \rightarrow
	H_{\operatorname*{loc}}^{1}\left(  \mathbb{R}^{n}\right),\\
	S&:H^{-1/2}%
	\left(  S_{R}\right)  \rightarrow H^{1/2}\left(  S_{R}\right)  ,\\
	T&:H^{1/2}\left(  S_{R}\right)  \rightarrow H^{1/2}\left(  S_{R}\right)  
\end{align*}
are continuous linear operators. 
In this way, well-posedness of
(\ref{eqn:exterior_problem}) with condition $\operatorname*{IC}\nolimits_{s}%
^{\operatorname*{Green}}$ follows if $S$ in (\ref{Seq}) has a bounded inverse.
Theorem 7.6 in \cite{Mclean00} implies that $S$ is a Fredholm operator of
index $0$ so that uniqueness implies bounded invertibility. Uniqueness follows
for $n\geq2$ and $s\in\CutC$ from
\cite[Theorem 9.11]{Mclean00}, for $n\geq3$ and $s=0$ from \cite[Corollary
8.11]{Mclean00} and for $n=2$, $s=0$, and $R\neq1$ from \cite[Theorem
8.16.(ii)]{Mclean00}. Hence (\ref{eqn:exterior_problem}) 
with condition $\operatorname*{IC}\nolimits_{s}%
^{\operatorname*{Green}}$ is well-posed in these three cases (with continuous 
solution operator).
This finishes the proof of \Cref{PropProblems}.1a and \ref{PropProblems}.2b, whereas \Part \ref{PropProblems}.1b
directly follows from the implication (\ref{GreenStrong}).

It remains to prove \Cref{PropProblems}.2a. Uniqueness
of solutions is guaranteed by
\cite[Thm~8.10]{Mclean00}. Existence and well-posedness follows 
from the expansion (\ref{uhatexpLap}) of the solution
\begin{equation}
\widehat{u}\left(  r,\xi\right)  =\sum_{m\in\mathbb{N}_{0}}\left(  \frac{R}%
{r}\right)  ^{m+2\nu}\sum_{j\in J_{m}}\widehat{g}_{m,j}Y_{m,j}\left(
\xi\right)  \quad\text{with\quad}u\left(  r\xi\right)  =\widehat{u}\left(
r,\xi\right)  \label{seriesuhat}%
\end{equation}
We prove that the mapping $g\rightarrow u$ induced by this series is a
continuous map from $H^{1/2}\left(  S_{R}\right)  $ to $H_{\operatorname*{loc}%
}^{1}\left(  B_{R}^{+}\right)  $. For this, we show that for any $d>R$ it holds%

\begin{equation}
\left\vert
\kern-.1em%
\left\vert
\kern-.1em%
\left\vert u\right\vert
\kern-.1em%
\right\vert
\kern-.1em%
\right\vert _{H^{1}\left(  B_{R}^{+}\cap B_{d}\right)  }:=\left(  \sum
_{m\in\mathbb{N}_{0}}\sum_{j\in J_{m}}\left(  F_{m,j}^{2}+G_{m,j}^{2}%
+H_{m,j}^{2}\right)  \left\vert \widehat{g}_{m,j}\right\vert ^{2}\right)
^{1/2}<\infty\label{wfourier}%
\end{equation}
with%
\begin{align*}
F_{m,j}^{2}  & :=\int_{R}^{d}\left(  \frac{R}{r}\right)  ^{2m+4\nu}%
r^{n-1}dr,\quad G_{m,j}^{2}:=\int_{R}^{d}\left\vert \frac{d}{dr}\left(
\frac{R}{r}\right)  ^{m+2\nu}\right\vert ^{2}r^{n-1}dr,\\
H_{m,j}^{2}  & :=\int_{R}^{d}\frac{m^{2}}{r^{2}}\left(  \frac{R}{r}\right)
^{2m+4\nu}r^{n-1}dr.
\end{align*}
These coefficients allow for the straightforward estimate:%

\begin{align*}
F_{m,j}^{2} &  \leq\left\{
\begin{array}
[c]{cc}%
\log\frac{d}{R} & m+\nu=0,\\
\frac{R^{2\nu}}{2(m+\nu)} & m+\nu\geq1,
\end{array}
\right.  ,\quad G_{m,j}^{2}\leq\left(  m+2\nu\right)  ^{2}\left\{
\begin{array}
[c]{ll}%
\frac{R^{2(\nu+1)}}{2(m+\nu-2)} & m+\nu\geq3,\\
R^{2\left(  \nu+1\right)  }\log\frac{d}{r} & m+\nu=2,\\
\frac{1}{2}\left(  d^{2}-R^{2}\right)  R^{2\nu} & m+\nu=1,\\
\frac{1}{4}\left(  d^{4}-R^{4}\right)  R^{2\nu-2} & m+\nu=0,
\end{array}
\right.  \\
H_{m,j}^{2} &  =\left\{
\begin{array}
[c]{ll}%
0 & m+\nu=0,\\
\frac{m^{2}R^{2\nu}}{2(m+\nu)} & m+\nu\geq1.
\end{array}
\right.
\end{align*}
Consequently,
\[
\max\left\{  F_{m,j}^{2},G_{m,j}^{2},H_{m,j}^{2}\right\}  \leq\left(  \frac
{d}{R}\right)  ^{4}R^{2\nu}\left(  1+R^{2}\right)  \left(  m+2\nu\right)
\left(  \nu+2\right)  .
\]
The combination of the previous arguments reveals
\begin{align*}
\left\vert
\kern-.1em%
\left\vert
\kern-.1em%
\left\vert u\right\vert
\kern-.1em%
\right\vert
\kern-.1em%
\right\vert _{H^{1}\left(  B_{R}^{+}\cap B_{d}\right)  } &  \leq\sqrt
{3}\left(  \frac{d}{R}\right)  ^{2}R^{\nu}\sqrt{\left(  1+R^{2}\right)
\left(  \nu+2\right)  }\left(  \sum_{m\in\mathbb{N}_{0}}\sum_{j\in J_{m}%
}\left(  m+2\nu\right)  \left\vert \widehat{g}_{m,j}\right\vert ^{2}\right)
^{1/2}\\
&  \leq\sqrt{3}\left(  \frac{d}{R}\right)  ^{2}R^{\nu}\sqrt{\left(
1+R^{2}\right)  \left(  \nu+2\right)  }\left(  1+2\nu\right)  \left\Vert
g\right\Vert _{H^{1/2}\left(  S_{R}\right)  }.
\end{align*}
This implies for any $g\in H^{1/2}\left(  S_{R}\right)  $ that the series
(\ref{seriesuhat}) defines a continuous map to a function $u\in
H_{\operatorname*{loc}}^{1}\left(  B_{R}^{+}\right)  $ and concludes the proof of
existence and well-posedness in \Cref{PropProblems}.2a.\qed

\medskip
\noindent\textbf{Data availability.} 
Not applicable. The manuscript has no associated data.%

\medskip
\noindent\textbf{Conflict of interest.} 
There is no conflict of interest.

\medskip
\noindent\textbf{Acknowledgements.}
The authors gratefully acknowledge the financial support by 
the Swiss National Science Foundation under grant no.\ 2000-1-240043.
\bibliographystyle{abbrv}
\bibliography{Bibliography,nlailu}

\newcommand{\noopsort}[1]{} \newcommand{\printfirst}[2]{#1} \newcommand{\singleletter}[1]{#1} \newcommand{\switchargs}[2]{#2#1} \def\cprime{$'$} \def\cprime{$'$} \def\cprime{$'$}
\begin{thebibliography}{10}

\bibitem{BLS:VariableOrderDirectional2021}
S.~B{\"o}rm, M.~{Lopez-Fernandez}, and S.~A. Sauter.
\newblock Variable order, directional {{H}}{$^2$}-matrices for {{Helmholtz}} problems with complex frequency.
\newblock {\em IMA J. Numer. Anal.}, 41(4):2896--2935, 2021.

\bibitem{MonkChandlerWilde}
S.~Chandler-Wilde and P.~Monk.
\newblock Wave-number-explicit bounds in time-harmonic scattering.
\newblock {\em SIAM J. Math. Anal.}, 39:1428--1455, 2008.

\bibitem{coltonkress_inverse}
D.~Colton and R.~Kress.
\newblock {\em Inverse {A}coustic and {E}lectromagnetic {S}cattering {T}heory}, volume~93 of {\em Applied Mathematical Sciences}.
\newblock Springer-Verlag, Berlin, 1992.

\bibitem{NIST:DLMF}
{\it NIST Digital Library of Mathematical Functions}.
\newblock http://dlmf.nist.gov/, Release 1.0.13 of 2016-09-16.
\newblock F.~W.~J. Olver, A.~B. {Olde Daalhuis}, D.~W. Lozier, B.~I. Schneider, R.~F. Boisvert, C.~W. Clark, B.~R. Miller and B.~V. Saunders, eds.

\bibitem{Fen:FiniteElementMethod1983}
K.~Feng.
\newblock Finite element method and natural boundary reduction.
\newblock In {\em Proc. {{Int}}. {{Congr}}. {{Math}}.}, pages 1439--1453, 1983.

\bibitem{Gilbarg}
D.~Gilbarg and N.~Trudinger.
\newblock {\em Elliptic {P}artial {D}ifferential {E}quations of {S}econd {O}rder}.
\newblock Springer-Verlag, 1983.

\bibitem{GPS:HelmholtzEquationHeterogeneous2019}
I.~Graham, O.~Pembery, and E.~Spence.
\newblock The {{Helmholtz}} equation in heterogeneous media: {{A}} priori bounds, well-posedness, and resonances.
\newblock {\em Journal of Differential Equations}, 266(6):2869--2923, Mar. 2019.

\bibitem{Gra:OptimalTraceNorms2025}
B.~Gr{\"a}{\ss}le.
\newblock Optimal trace norms for {{Helmholtz}} problems.
\newblock {\em arXiv:2506.11944}, 2025.

\bibitem{GHS:StableSkeletonIntegral2025}
B.~Gr{\"a}{\ss}le, R.~Hiptmair, and S.~A. Sauter.
\newblock Stable skeleton integral equations for general coefficient {{Helmholtz}} transmission problems.
\newblock {\em arXiv:2507.00991}, 2025.

\bibitem{KG:ExactNonreflectingBoundary1989}
J.~B. Keller and D.~Givoli.
\newblock Exact non-reflecting boundary conditions.
\newblock {\em Journal of Computational Physics}, 82(1):172--192, May 1989.

\bibitem{LionsMagenesI}
J.~Lions and E.~Magenes.
\newblock {\em Non-{H}omogeneous {B}oundary {V}alue {P}roblems and {A}pplications}.
\newblock Springer-Verlag, Berlin, 1972.

\bibitem{MM:FiniteElementMethod1980}
R.~C. Maccamy and S.~P. Marin.
\newblock A finite element method for exterior interface problems.
\newblock {\em International Journal of Mathematics and Mathematical Sciences}, 3(2):311--350, Jan. 1980.

\bibitem{Mclean00}
W.~McLean.
\newblock {\em Strongly Elliptic Systems and Boundary Integral Equations}.
\newblock Cambridge, Univ. Press, 2000.

\bibitem{MelenkSauterMathComp}
J.~M. Melenk and S.~A. Sauter.
\newblock Convergence analysis for finite element discretizations of the {H}elmholtz equation with {D}irichlet-to-{N}eumann boundary conditions.
\newblock {\em Math. Comp}, 79:1871--1914, 2010.

\bibitem{Nedelec01}
J.~C. N{\'e}d{\'e}lec.
\newblock {\em Acoustic and {E}lectromagnetic {E}quations}.
\newblock Springer, New York, 2001.

\bibitem{Sch:EightyYearsSommerfelds1992}
S.~H. Schot.
\newblock Eighty years of {{Sommerfeld}}'s radiation condition.
\newblock {\em Historia Mathematica}, 19(4):385--401, Nov. 1992.

\bibitem{Shubin2001}
M.~A. Shubin.
\newblock {\em Pseudodifferential {O}perators and {S}pectral {T}heory}.
\newblock Springer-Verlag, Berlin, second edition, 2001.

\bibitem{SW:WavenumberExplicitParametricHolomorphy2023}
E.~A. Spence and J.~Wunsch.
\newblock Wavenumber-explicit parametric holomorphy of {{Helmholtz}} solutions in the context of uncertainty quantification.
\newblock {\em SIAM/ASA J. Uncertainty Quantification}, 11(2):567--590, June 2023.

\bibitem{Watson}
G.~N. Watson.
\newblock {\em A Treatise on the Theory of Bessel Functions}.
\newblock Cambridge University Press, 1966.

\end{thebibliography}
\end{document}